\documentclass[11pt]{amsart}
\usepackage[english]{babel}
\usepackage[colorlinks,linkcolor=red,anchorcolor=green,citecolor=blue]{hyperref}
\usepackage{booktabs}
\usepackage{array} % for better arrays (eg matrices) in maths
\usepackage{paralist} % very flexible & customisable lists (eg. enumerate/itemize, etc.)
\usepackage{verbatim} % adds environment for commenting out blocks of text & for better
\usepackage{subfig} % make it p
\usepackage{bm}
\usepackage{graphicx} % support the \includegraphics command and options
\usepackage{amsmath,breqn}
\usepackage{geometry} % to change the page dimensions
\usepackage{amssymb}
\usepackage{amsthm}
\usepackage{float}
\usepackage{amsfonts,amssymb}
\usepackage{color, soul}
\usepackage{makecell}
\usepackage{colortbl,dcolumn}
\usepackage{diagbox}
\usepackage{stfloats}
\usepackage{enumerate}
\usepackage{epstopdf}
\usepackage{tabularx}
\usepackage{caption}
\usepackage{siunitx}
\usepackage{multirow}
\usepackage{xmpmulti}
\usepackage{arydshln}
\theoremstyle{definition}

\def\eop{\hfill\square}

\newtheorem{Def}{Definition}
\newtheorem{Thm}{Theorem}

\newtheorem{Exmp}{Example}
\title[]{Sparse Representations of Solutions to a class of Random Boundary Value Problems}
\author[Fang Yang]{Fang Yang}
\address{Fang Yang, Department of Mathematics\\
Zhejiang Normal University\\
China}
\email{19957019632@189.cn}
\author[Min Chen]{Min Chen}
\address{Min Chen, School of Mathematical Sciences, Shanxi University.
\newline
Academy of Mathematics and Systems Science, Chinese Academy of Sciences\\
China}
\email{mchen@amss.ac.cn}
\author[Jiecheng Chen]{Jiecheng Chen}
\address{Jiecheng Chen, Department of Mathematics\\
Zhejiang Normal University\\
China}
\email{jcchen@zjnu.edu.cn}
\author[Ying Zhang]{Ying Zhang}
\address{Ying Zhang, Macau Institute of Systems Engineering\\
Macau University of Science and Technology\\
China}
\email{cnuzhangying@163.com}

\author[Tao Qian]{Tao Qian*}
\address{Tao QIAN, Macao Center for Mathematical Sciences\\
Macau University of Science and Technology\\
China}
\email{tqian@must.edu.mo}
\thanks{*Corresponding author.}

\begin{document}
\maketitle
\begin{abstract}
We introduce certain sparse representation methods, named as stochastic pre-orthogonal adaptive Fourier decomposition 1 and 2 (SPOAFD1 and SPOAFD2) to solve the Dirichlet boundary value problem and the Cauchy initial value problem of random data.  To solve the stochastic boundary value problems the sparse representation is, as the initial step, applied to the random boundary data. Due to the semigroup property of the Poisson and the heat kernel, each entry of the expanding series can be lifted up to compose a solution of the Dirichlet and the Cauchy initial value problem, respectively. The sparse representation gives rise to analytic as well as numerical solutions to the problems with high efficiency.
\end{abstract}

\bigskip


\bigskip

\date{today}
%\tableofcontents

\section{Introduction}

In many science and engineering problems, the stochastic partial differential equations (SPDEs)
%and the partial differential equations (PDEs) with the random coefficient/data
 are more realistic than deterministic PDEs. It is particularly important to give effective numerical solutions of SPDEs (\cite{Kloeden-1992,Milstein-1995,Milstein-2004,Maruyama-Euler-Method-1955,Allen_1998, Debussche_2007,Walsh_Finite-Element-2005,Machiels_Fourier-spectral-1998,
 Hausenblas_Galerkin-method2003}).
It is assumed that $(\Omega,\mathcal{F}, \mathbb{P})$ is a probability space and $D$ is a region of $\mathbb{R}^n$ throughout this paper. In this paper we restrict ourselves to deal with the cases where $D$ is either the unit ball or the upper-half space. We will concern elliptic PDEs with random coefficients or random forcing terms
(\cite{LPS, Xiu2010-numerical-methods}), which have the form
\begin{equation}
\label{random-elliptic-Equ}
\left\{\begin{aligned}
&\nabla(\kappa(x,\omega) \nabla u(x,\omega) )=f(x,\omega),~ x\in D , {\rm a.s.}~\omega, \\&
u(x,\omega)=g(x,\omega), ~ x\in \partial D,~ {\rm a.s.}~\omega,
\end{aligned}\right.
\end{equation}
where
$  \kappa, f$ and $ g$ are $\mathbb{R}$-valued
functions on $D \times \Omega$.
For such SPDEs a number of numerical methods  have been developed, by using methods including the Monte Carlo finite element method, Karhunen-Lo\`eve (K-L) expansions, finite element methods (FEMs), stochastic Galerkin FEM and stochastic collocation FEM  (\cite{LPS,Xiu_2002, Mugler_2013,Babuka_Stochastic-Collocation2007}), etc.

Besides the above type SPDEs the method we will introduce may also be used to the parabolic type PDEs with random boundary- or/and initial-values, including, in particular, incompressible Navier Stokes equations (\cite{Sapsis_DO-2009}). They have the form
\begin{equation}\label{parabolic-SPDE}
\left\{
\begin{aligned}
&\partial_t u(t,x,\omega)=\mathcal{L}[ u(t,x,\omega),\omega], ~ ~x \in D , ~~ t \in \mathbb{R}_+,~{\rm a.s.}~\omega \in \Omega,\\
&u(0,x,\omega)=u_0(x,\omega), ~ ~x \in D ,~{\rm a.s.}~\omega \in \Omega, \\
&\mathcal{B}[u(t,x,\omega)]=h(t,x,\omega),~x \in \partial D,~{\rm a.s.}~\omega \in \Omega,
\end{aligned}
\right.
\end{equation}
where $ \mathcal{L} $ is a linear or non-linear differential operator, $ x \in D $ is the spatial coordinate, $ t $ is the time variable, and $ \mathcal{B} $ is a
linear differential operator.
%For the ease of notation in what follows we write the explicit dependence of $ \mathcal{L} $ on $ (t,x,\omega) $.
In the classical method the problem (\ref{parabolic-SPDE}) can be reduced to a certain type of deterministic PDEs and a system of stochastic ordinary differential equations by the so called dynamically orthogonal representation proposed in \cite{Sapsis_DO-2009,Choi_DO-schemes-2013}.

Solutions of SPDEs are again random processes. Hence, it is meaningful to study how to decompose a random process explicitly with efficient convergence. Before we introduce stochastic adaptive Fourier decomposition (SAFD) type sparse representations (\cite{QianSAFD}), we first recall two classical techniques for representing random processes (\cite{Nouy_Recent-developments-2009}).

Among Wiener chaos spectral expansions there is Fourier-Hermite expansion being a classical technique for Brownian motion  (\cite{Hou_2006,Yannacopoulos_2011,Kaligotla_2012,Jardak_2002,Luo_2006}). Assume $ \xi(\omega)\triangleq\{\xi_i(\omega)\}_{i=1}^n $ is a set of independent standard Gaussian random variables. For a finite index
$  \alpha =(\alpha_1,\alpha_2,\cdots,\alpha_n)\in \mathbb N^n $
, the multivariate Hermite polynomial of $x=(x_1,\cdots,x_n)\in \mathbb R^n$ can be written as:
\begin{equation}\label{Hermite-poly-n-D}
H_\alpha(x)\triangleq\prod^{n}_{i=1}H_{\alpha_i}(x_i),
\end{equation}
where each $H_{\alpha_i}(x_i)$ is a one dimensional Hermite polynomials of $x_i.$ One dimensional Hermite polynomials form an orthonormal system
\[ \langle h_n,h_m \rangle_{\nu}
\triangleq\int_{-\infty}^{+\infty} h_n h_m d\nu=\delta^n_m \]
that induces the orthonormality of the higher dimensional systems.
Suppose that $u(\xi(\omega))$ is a function of the random
variables $\xi(\omega)$ with $\mathbb{E}_\omega|u(\xi(\omega))|^2<+\infty$, then the function  $ x\mapsto u(x)$ belongs to $L^2(\mathbb{R}^n, d\nu^n)$ and
$u(\xi(\omega))$ having a Fourier-Hermite expansion converging in $L^2(\Omega , d\mathbb{P})$:
\begin{equation}\label{Wiener-chaos-expansion-1D}
u(\xi(\omega))=\sum_{\alpha }u_\alpha
H_{\alpha}(\xi(\omega)),
\end{equation}
where the Fourier coefficients  $u_\alpha\triangleq\mathbb{E}[ u(\xi(\cdot))  H_{\alpha}(\xi(\cdot))].$

For $x \in D$, assume $ u(x,s,\omega)\triangleq f(x, W_s(\omega))$ is a function of the Brownian motion $\{W_s(\omega)\}_{0\leq s\leq T}$ with zero means. Analogous to the finite dimensional Fourier-Hermite expansion (\ref{Wiener-chaos-expansion-1D}), it admits a similar expansion. Before we state it
we recall the definition of the second order processes:
\begin{Def}\label{second}(\cite{LPS})
	A stochastic processes $\{u(x,\omega)\}_{x \in D}$ is called a second order  processes if for almost everywhere $ x \in D  $,~
	\begin{equation}\label{def-second-order-prcs}
	\mathbb{E}|u(x,\cdot)|^2=\int_{\Omega}|u(x,\omega)|^2 d\mathbb{P}(\omega)
	< \infty.
	\end{equation}
\end{Def}
\noindent
Define the set of finite multi-indices as
\[
\mathcal{J} \triangleq\{
\alpha\in {\mathbb N}^{{\mathbb N}_+} :~|\alpha|=\sum_i \alpha_i<\infty
\}.
\]
For $\alpha \in \mathcal{J}$, define the Hermite polynomial  of $\xi(\omega)$ similar to (\ref{Hermite-poly-n-D}) by the tensor product
\begin{equation}
\label{Wick-poly}
T_\alpha(\xi(\omega))\triangleq\prod_{i=1}^{\infty}H_{\alpha_i}(\xi_i(\omega)),
\end{equation}
which is called a random Wick polynomial.
The above product comprises only finite number of factors since $|\alpha|<\infty$. Given a second order process $ u(x,s, \omega)\triangleq f(x,W_s(\omega)), $ assumed to be a function of the Brownian motion $\{W_s(\omega)\}_{0\leq s\leq T}$, Cameron and Martin (\cite{Luo_2006, Cameron-Martin1947}) provided a Fourier expansion as the following,
\begin{equation}\label{WCE}
u(x,s,\omega)=\sum_{\alpha \in \mathcal{J} }u_\alpha(x,s) T_\alpha(\xi(\omega)),
\end{equation}
in $L^2(\Omega, d\mathbb P)$,
where $T_{\alpha}(\xi(\cdot))$ is defined by (\ref{Wick-poly}) and  $u_\alpha(x,s)\triangleq\mathbb{E}[ u(x,s,\cdot)  T_{\alpha}(\xi(\cdot))]$.
The above expansion is usually called Wiener chaos expansion of $ u(x, s,\omega)$ that represents the randomness of the solution analytically by a set of random bases with deterministic coefficients. This method reduces the SPDEs to solving the deterministic PDEs of the coefficients (\cite{Luo_2006, Georgi_2017}).

The other common approach for expanding random processes is K-L expansion (\cite{LPS}), which applies to second order processes $ u$ restricted to a space $L^2(\Omega, L^2(D)).$
\begin{Def} (\cite{LPS,QianSAFD})
	Suppose $D \subset \mathbb R ^n.$ The space $ L^2(\Omega, L^2(D)) $ is defined to be the Hilbert space consisting of all $ L^2(D) $-valued random variables $v: D \times \Omega \rightarrow \mathbb{R}$ satisfying
	\begin{equation}
	\|v\|_{L^2(\Omega, L^2(D))}^2\triangleq
	\int_{\Omega} \int_{D} |v(x, \omega)|^2 dx d\mathbb{P}<\infty.
	\end{equation}
\end{Def}
\def\N{\mathcal N}
\noindent
For brevity, we also write $ \mathcal{N}=L^2(\Omega,L^2(D)).$ The reference \cite{QianSAFD} uses the same type of spaces on the classical domains defined in terms of Fourier expansions with random coefficients, being equivalent with the above definition due to the Plancherel Theorem. \cite{QianSAFD} develops stochastic sparse kernel-representations of stochastic reproducing kernel Hilbert spaces and stochastic Hilbert spaces with a dictionary.  The spaces $ L^2(\Omega, L^2(D))$ have also been studied as Bochner spaces which are Hilbert spaces with an inner product defined through Bochner integrals. See \cite{Hytonen_AnalysisBanachSpaces_2016, Yosida1978_FunctionalAnalysis, Schwab2003_sparse-FE-PDE-stochastic} and the references therein.
%\begin{equation}
%\label{Canonical-iso}
%L^2(\Omega, \mathcal H(D))\cong
%\mathcal H(D)\otimes L^2(\Omega, d\mathbb{P})\cong
%\mathcal H(D,L^2(\Omega, d\mathbb{P})),
%\end{equation}
%where $\mathcal H(D)\otimes L^2(\Omega, d\mathbb{P})$ denotes the tensor product of $\mathcal H(D)$ and $ L^2(\Omega, d\mathbb{P})$.

 Let $u(x,\omega)$ be a process in  $L^2(\Omega, L^2(D)) $ and $ C(x,y)\triangleq{\rm Cov}(u(x,\cdot),u(y,\cdot))$ the covariance of the $u.$
 By using $C(x,y)$ one defines the integral operator $ \mathcal{ C } $
\begin{equation}\label{cov-operator}
\mathcal{ C }(\phi)(x)\triangleq
\int_D C(x,y)\phi(y)dy,~x\in D,~\phi \in L^2(D).
\end{equation}
Setting $\tilde{u}(x)\triangleq\mathbb{E}(u(x,\cdot)),$ we have the K-L expansion in the sense of $ L^2(\Omega, L^2(D))$,
\begin{equation}\label{K-L-expansion-L2}
u(x,\omega)=\tilde{u}(x)+\sum_{j=1}^{\infty}\sqrt{v_j} \phi_j(x)\xi_j(\omega),
\end{equation}
where the Fourier coefficients
$ \xi_j(\omega)=\frac{1}{\sqrt{v_j}}
\langle {u(\cdot,\omega)-\tilde{u}(\cdot),\phi_j(\cdot)}
\rangle_{L^2(D)} $ ,
and
$\{v_j, \phi_j\}_j$ are the eigen-pairs
of the operator $\mathcal{ C }$. Combined with the Galerkin method, the K-L Galerkin procedure is a type of Galerkin method that employs, as basis functions, the empirical eigenfunctions of the K-L decomposition (\cite{Park1998_K-L, Park2000_K-L}). It consequently reduces the original PDE to a small set of ordinary differential equations or algebraic equations. A boundary optimal control problem of the Navier-Stokes equation is solved by this technique in \cite{Park2000_K-L}. Several complexity reduction methods have been proposed, such as model reduction methods for SPDEs (\cite{Boyaval_2009,Nouy_spectral-decpst2007,Nouy_spectral-decpst2008}), sparse approximation
methods for high-dimensional approximations (\cite{Cohen-2010,Schwab-2011,Nobile-2008,Todor-2007}).

Both the above reviewed Fourier-Hermite expansion and the Karhunen-Lo\`eve decomposition offer series expansions of the considered random signal or random process into entries being of the tensor product form. Based on the tensor forms one is able to solve the stochastic boundary value or initial value problems with deterministic mathematical analysis methodology. The decompositions are just general and not specially made for solving particular SPDE problems. They, therefore, do not offer convenience to particular questions, and exhibits low efficiency. In the present study, we propose an easily implementable sparse representation method, based on optimal selections of dictionary elements made of fundamental solutions of the problems. This, through the $\lq\lq$lifting up" method based on the semigroup property of the fundamental solutions (see \S 4), offers remarkable convenience and efficiency to solve the SPDE problems. In many cases, the dictionary elements are parameterized reproducing kernels of a relevant analytic functional Hilbert spaces.

 The methodology of this paper is based on the theory and algorithms of Qian in \cite{QianSAFD} in which SAFD for the Hardy space is first developed and then further extended to become stochastic pre-orthogonal adaptive Fourier decomposition 1 and 2 (SPOAFD1, SPOAFD2) for general Hilbert spaces with a dictionary.  The SAFD type methods are further development of what is called adaptive Fourier decomposition (AFD) whose comprehensive study includes a number of variations including  core adaptive Fourier decomposition (or just AFD, see
\cite{QWa,Qianbook}), pre-orthogonal adaptive Fourier decomposition (POAFD, available for general Hilbert spaces with a dictionary in which a Takenaka-Malmquist system is unavailable, see \cite{QSW,Q2D,qu2018,qu2019}), unwinding Blaschke expansion (unwinding AFD, independently developed also by Coifman et al., see \cite{CS,CP,Qian2010,QTa}, $n$-best reproducing kernel approximation ($n$-best AFD, see \cite{QWM,QQ}), and also in higher dimensions for matrix-valued functions (\cite{ACQS1,ACQS2}). The study of the AFD type algorithms was originated from signal decomposition into meaningful positive instantaneous frequency in physics (\cite{QWa}). The maximal selection principle of the AFD methods shares the matching pursuit idea (or greedy algorithms) but emphasises attainability of the global optimal parameters that leads to the delicate analysis on multiple selection of parameter and the multiple kernel concept. The contribution of the present study is to initiate application of the SAFD methods that gives solutions of the stochastic boundary value problems in the form of fast converging series of lifted and shifted fundamental solutions of the related operators. \\

%摘自Nouy_Random fields representations for stochastic elliptic boundary value problems and statistical inverse problems2013

%%%%%%%%%%%%%%%%%%%%
We note that the present paper deals with two types of random data on the boundary. The first type is those expressible as a bi-variate-function, $f(t,X),$ of which one is the boundary variable $t$ and the other is a random variable $X$. The second type is that the random data on the boundary is actually a random process. For the latter we take Brownian bridge $W_t$ as an example. Detailed studies involving comparison of various types of decomposition methods and applications of the SAFD methods are reported in our separate papers (\cite{KLvsAFD},\cite{OUprocess}).\\

Denote by $B_r(a)$ the open ball centred at $a$ with radius $r$ in $\mathbb{R}^{n+1},$ and $\partial B_r(a)$ its boundary, also called the $r$-sphere.  In particular, we write $B_1\triangleq B_1(0)$ in short.
We deal with the Dirichlet problem for Laplace equation with a random boundary condition
%%%%我们考虑的问题分别是椭圆型、抛物型pde中最简单的，但是与问题（3）的不同点在于，我们让初边值是二阶过程（随机信号）。
\begin{equation}\label{Laplace-equ-random}
\left\{
\begin{aligned}
&\Delta_x u(x,\omega)=0, ~ ~x \in B_1\subseteq \mathbb{R}^{n+1},~{\rm a.s.}~\omega \in \Omega,\\
&u(x,\omega)=f(x,\omega),~ ~x\in{\partial  B_1},~{\rm a.s.}~\omega \in \Omega,
\end{aligned}
\right.
\end{equation}
and the Cauchy initial problem of the heat equation with a random initial condition
\begin{equation}\label{heat-equ-random}
\left\{
\begin{aligned}
&(\partial_t-\Delta_x )u(t,x,\omega) =0, ~ ~(t,x) \in \mathbb{R}^{1}_+\times \mathbb{R}^n,~{\rm a.s.}~\omega \in \Omega,\\
&u(0,x,\omega)=f(x,\omega),~ ~x \in \mathbb{R}^n,~{\rm a.s.}~\omega \in \Omega.
\end{aligned}
\right.
\end{equation}
%%%%%%%%%%%
Becus and Cozzarelli defined the random generalized solutions of such types of SPDE problems, and proved existence and uniqueness of the solutions (\cite{Becus_rdm-Laplace-1976,Becus_rdm-Laplace-1977,Becus_rdm-heat-1977}).
Becus also proposed an iterative scheme for solving such random heat equation and established its convergence (\cite{Becus_approximations-rdm-heat-Equ-1978}). Tasaka proposed the finite element solutions of the one-dimensional heat equation with random initial conditions (\cite{Tasaka_FE-method-Rdm-heat-Equ-1983}).

In this paper, using the random signal analysis methodology theoretically developed in \cite{QianSAFD}, we give and analyse the solutions of
(\ref{Laplace-equ-random}) and (\ref{heat-equ-random}).
Different from the existing methods, the proposed SAFD methods are easy to be implemented and directly give solutions of the SPDE problems under study for all random boundary or initial data. Our methods depend on the distribution of $X$ for the boundary bi-variate-function $f(t,X)$ type, and the covariance function $C(s,t)$ for random process type. For the latter our algorithm does not require any information on the related eigenvalues, nor the eigenfunctions, as it does for proceeding the K-L decomposition.\\

The paper is organized as follows. In \S \ref{Preliminaries} we review relevant preliminary knowledge and define notations that will be used in the following sections. In \S \ref{Sect-random-PDEs} we formulate the Dirichlet problem of Laplace equation with a random boundary condition and the Cauchy problem of heat equation with a random initial condition, and, respectively, introduce the stochastic harmonic Hardy space and the stochastic heat Hardy space as the appropriate function spaces for the solutions. In \S \ref{Sect-SPOAFD} we develop the SPOAFD methodology in obtaining sparse representations of solutions of the concerned stochastic boundary/initial problems. The stochastic boundary vanishing conditions for the given function spaces are proved, which guarantee applicability of the SPOAFD algorithms.
The numerical experiments are presented in \S \ref{Sect-experiment} to confirm the theoretical estimates.

%%%对于SPDEs中的白噪声项，也有两种近似方法。第一种是Wiener展开[], 第二种是正则化逼近
\def\D{D}
\section{Preliminaries}
\label{Preliminaries}
 We will first set our notation for the deterministic case. For the Dirichlet problem of Laplace equation we focus on the $(n+1)$-dimensional unit ball $B_1(0)$ and the upper half-space $\mathbb{R}^{1+n}_{+}$ contexts,  and for the Cauchy initial problem of the heat equation on $ \mathbb{R}_+^{1+n} $.

For deterministic Dirichlet problems in bounded domains we adopt the formulation that $\D\subset \mathbb  R^{1+n}$ is a bounded domain with a smooth boundary, and $f \in L^2(\partial \D ).$
We need to determine a function $u$ satisfying
\begin{equation}\label{Prob-Dirichlet}
\left\{
\begin{aligned}
&\Delta u(x)=0, ~ \forall~x \in \D,\\
& u(x)_{\partial \D}=f(x), \ {\rm for\ a.e.}\ x \in{\partial \D}.
\end{aligned}
\right.
\end{equation}
We note that the second equation in (\ref{Prob-Dirichlet}) takes the non-tangential boundary limit sense. The solution to (\ref{Prob-Dirichlet}) can
be written as the boundary integral (\cite{Han-Lin_PDEs2011})
\begin{equation}\label{Green-funct-reprst}
u(x)=\int_{\partial \D} \frac{\partial }{\partial n_y} G(x,y) f(y) d\sigma(y),
\end{equation}
where $G(x,y)$ is the Green's function of $D$. Restricting the question to the unit ball $B_1= \D$ %or the half-space $\mathbb{R}^{1+n}_{+}=\{(t,x)\in {R}^{1+n}:t>0 \}$%
, the function $  \frac{\partial }{\partial n_y} G(x,y)  $ becomes the Poisson kernel (for the unit ball) $P_x(y)$ defined as:
 %or $\mathbb{R}^{1+n}_{+}$

\begin{equation}\label{Poiss-Ker-B_1}
P_x(y^{\prime})\triangleq P(x,y^{\prime})\triangleq c_n \frac{1-|x|^2}{|x-y^{\prime}|^n},~
x \in B_1,~ y^{\prime} \in \partial B_1,
\end{equation}
%or
%\begin{equation}
%P(x,y):=c_n \frac{1-|x|^2}{|x-y|^n}
%\end{equation}
where $c_n$ is a constant depending on $n$ that can be different from line to line.
Thus the solution of the problem can be written explicitly, as Poisson integral of $f,$
\begin{equation}
u_f(x)\triangleq\int_{\partial B_1} P(x,y^{\prime}) f(y^{\prime}) d\sigma(y^{\prime}).
\end{equation}
We can also obtain the explicit formula for the solution of (\ref{Prob-Dirichlet})  for the upper half-space $\mathbb{R}^{1+n}_+=\D$ and $f\in L^2(\mathbb{R}^n),$
\begin{equation}
u_f(x)=\int_{\mathbb{R}^n} {P} (x,\underline{y}) f(\underline{y}) d\underline{y},\end{equation}
where ${P} (x,\underline{y})$ is the counterpart Poisson kernel for the upper half-space $ \mathbb{R}_+^{1+n},$ where $x$ is often written as $x=x_0+\underline{x}, x_0\in \mathbb{R}^{1}_+, ~\underline{x}\in \mathbb{R}^{n},$ and
\[{P} (x,\underline{y})=
c_n\frac{x_0}{(x_0^2+|\underline{x}-\underline{y}|^2)^{\frac{n+1}{2}}}.\]

For
$\D$ being either the unit ball $B_1$ or the upper-half space case there exists a harmonic Hardy $h^p$ space theory.
%The following results give the $L^p$ estimate of the Poisson integral when the boundary data $f \in L^p(\partial B_1)$.
We give an account for the unit ball case. For the purpose of this study, we restrict to $p=2.$
 The harmonic Hardy space  $h^2(B_1)$ is defined to be
\begin{Def} (\cite{Axler_Harmonic-funct-thrm})
(The harmonic Hardy space on the ball)
	\begin{equation}
	h^2(B_1)\triangleq\{u: \Delta u|_{B_1}=0,~\sup_{0<r<1}\| u(r\cdot) \|_2<\infty   \}.
	\end{equation}
The norm of $u\in h^2(B_1)$ is defined as $\|u\|_{h^2}=\sup\limits_{0<r<1}\| u(r\cdot) \|_2.$
\end{Def} \noindent

It is a fundamental result that every $u \in h^2(B_1)$ has a non-tangential boundary limit function (\cite{Stein-Weiss1971}) and there holds
\[ \sup_{0<r<1}\| u(r\cdot) \|_2=\lim_{r\rightarrow 1}\|u(r\cdot) \|_2
=\|\lim_{r\rightarrow 1}u(r\cdot) \|_2
, ~\forall~u \in h^2(B_1).
\]
The $h^2(B_1)$ is a Hilbert space equipped with the inner product
\begin{equation}\label{product-h2}
\langle u,v\rangle\triangleq\int_{\partial B_1}f(s)g(s)d\sigma(s),~\forall~ u,v\in h^2(B_1),
\end{equation}
where the $f(s),~g(s)$ are the non-tangential limits of, respectively, the Poisson integrals $u=P[f],~v=P[g],$ and $d\sigma(s)$ is the normalized measure on the sphere. Thus,  $h^2(B_1)$
is isometric to $L^2(\partial B_1)$ , denoted
\begin{equation}\label{Isometric-L2-h2}
L^2(\partial B_1) \cong h^2(B_1).
\end{equation}
The harmonic Hardy space $h^2(B_1)$ is a special case of the Hardy space introduced by Stein and Fefferman (\cite{Fefferman-Stein_Hardy-spc-1972}). Within the same category, another Hardy space, $H^2_{heat},$  is defined in relation to the heat equation in $\mathbb{R}_{+}^{1+n}.$
Denote by
 $$ \varphi_t(x)\triangleq(\frac{1}{4t\pi})^{n/2}e^{-\frac{1}{4t}|x|^2}, $$
the heat kernel. Then a solution of the corresponding heat equation can be written as
\begin{eqnarray}\label{heat}  u_f(t,x)
\triangleq\int_{\mathbb{R}^n} \varphi_t(x-y)f(y)dy=(\varphi_t \ast f) (x) ,~f \in L^2(\mathbb{R}^n),
 \end{eqnarray}
and $u=u_f$ belongs to heat-Hardy space $ H^2_{heat}(\mathbb{R}_{+}^{1+n})$ defined as
 \begin{equation}\label{heart-Hardy-spc}
\{u:~\exists~f \in L^2(\mathbb{R}^n),~u=\varphi_t \ast f,~
\varphi^+(f)\triangleq
 \sup_{t>0}|(\varphi_t \ast f)|
\in L^2(\mathbb{R}^n)
\},
 \end{equation}
 with the norm $\| u\|_{H^2}\triangleq\| \varphi^+(f)\|_{L^2}\simeq \|f\|_{L^2}$.
Similarly,  $H^2_{heat}(\mathbb{R}_{+}^{1+n})$
is isometric to $L^2(\mathbb{R}^n)$. The point of such formulation is that the spaces $h^2(B_1)$ and $H^2_{heat}(\mathbb{R}_{+}^{1+n})$ are reproducing kernel Hilbert spaces.

We recall the needed reproducing kernel structure applicable to both $h^2(B_1)$
and $H^2_{heat}(\mathbb{R}_{+}^{1+n}).$

\begin{Def} (\cite{Saitoh_RKHS-2016}) (Reproducing kernel Hilbert space)
	Let $D\neq \emptyset$.
	A reproducing kernel Hilbert space (RKHS) on the set $D$ is a Hilbert space $\mathcal{H}(D) $ with a function
	\[ K: D \times D\rightarrow \mathbb{C} \]
	possessing the reproducing property	
\begin{equation}\label{reproducing-prop}
\left\{
\begin{aligned}
	&K_q\triangleq K(\cdot,q)\in \mathcal{H}(D),~\forall~ q \in D,\\
	&\langle f,K_q \rangle=f(q),~\forall~ q \in D,~\forall~ f \in \mathcal{H}(D).
\end{aligned}\right.
\end{equation}
The function $K_q(\cdot)$ in (\ref{reproducing-prop}) is called the reproducing kernel of $\mathcal{H}(D) $ at $q$.
\end{Def}

For instance, denote by $P_{t}(s^{\prime})$ the Poisson kernel defined by (\ref{Poiss-Ker-B_1}), where $ t\triangleq|t|t^{\prime},~t^{\prime} \in D=B_1,$ and $s\triangleq|s|s^{\prime},~  s^{\prime}\in B_1.$ The harmonic Hardy space $h^2(B_1)$ is a RKHS (\cite{Qu-Sps-represt-dirac}) with the  reproducing kernel
\begin{equation}\label{rpdc-ker-h2}
K(s,t)\triangleq\langle
P_s,
P_t
\rangle
=P_{|s||t|s^{\prime}}(t^{\prime})
=P_{|s||t|t^{\prime}}(s^{\prime}).
\end{equation}
The reproducing kernel for the heat equation is
\begin{equation}\label{rpdc-ker-h2}
K(x,y)\triangleq\langle
h_x,h_y
\rangle
=h_{s+x}(\underline{y})
=h_{t+y}(\underline{x}),
\end{equation}
where $x=t+\underline{x},~ y=s+\underline{y}$ both are in $D=\mathbb R^{n+1}_+,$ while $t,s$ are in $\mathbb R^{1}_+.$

\section{Dirichlet problem with stochastic boundary value and Cauchy problem with stochastic initial value}\label{Sect-random-PDEs}
Consider a Dirichlet problem with a random boundary data, i.e.
\begin{equation}\label{Equ-1}
\left\{
\begin{aligned}
&\Delta u(x,\omega)=0, ~ \forall~x \in D\subseteq \mathbb{R}^{n+1},~{\rm a.s.}~\omega \in \Omega,\\
&u(x,\omega)=f(x,\omega),~ \ {\rm for \ a.e.}\ ~x\in{\partial  D},~{\rm a.s.}~\omega \in \Omega,
\end{aligned}
\right.
\end{equation}
where $f \in  L^2(\Omega,  L^2({\D})),~ {\D}=B_1.$
We also  write $f_\omega(\cdot)$  as $f(\cdot,\omega)$ when no ambiguity arises.
Babuska was among the first to study rigorously existence of solutions of  Dirichlet problems with random boundary condition (\cite{Babuska_1961}). Becus and Cozzarelli studied existence and properties of general solutions of (\ref{Equ-1}), see \cite{Becus_rdm-Laplace-1976,Becus_rdm-Laplace-1977,Becus_rdm-heat-1977}.
Furthermore, explicit formulae for the correlation tensor of the generalized solutions are obtained when the random boundary function is a white noise or a homogeneous random field on a sphere (\cite{Sabelfeld2008-expansion-random-boundary-PDE}).

In fact, the equation (\ref{Equ-1}) can be regarded as a family of equations labeled with $\omega \in \Omega$. For each $\omega \in \Omega$ possibly excluding an event of zero probability,
\begin{equation}\label{Equ-2}
\left\{
\begin{aligned}
&\Delta u(x,\omega)=0, ~ ~x \in B_1,\\
& u(x)=f_\omega(x),~ ~{\rm a.e.}~x \in \partial B_1,
\end{aligned}
\right.
\end{equation}
where  $f \in  L^2(\Omega,  L^2(D)).$
Hence (\ref{Equ-2}) becomes the classical Dirichlet boundary problem. Thus we can consider a stochastic process formed by the classical \textit{Poisson integral}
\begin{equation}\label{solut-Poiss-Integl}
u_{f_\omega}(x)\triangleq \int_{\partial B_1}p(x,t^{\prime})f_\omega(t^{\prime})d\sigma(t^{\prime}),~x \in B_1,
\end{equation}
where $p(x,t^{\prime})$ is the Poisson kernel in $B_1$.
For a.s. $\omega \in \Omega$, it follows that $u_{f_\omega}(x)$
solves the (\ref{Equ-2}),~and
\[ \lim_{\rho\rightarrow1}u_{f_\omega}(\rho t^{\prime} )=f_\omega(t^{\prime} ) \]
in both the $ L^2$-norm and the a.e. pointwise sense.

The stochastic process $u_{f_\omega}(x)$ defined in (\ref{solut-Poiss-Integl}) solves the problem (\ref{Equ-1}) in the $\mathbb{P}-$ a.s. sense.~The solution of (\ref{Equ-2})  is  a.s.  in $h^2(B_1)$ and, space-wise,  it is isometrically  identical with the boundary value $f_{\omega} \in L^2(\partial B_1).$

Similarly, the Cauchy problem with a
stochastic initial value is proposed in $\mathbb{R}_+^{n+1},$
\begin{equation}\label{stoc-heat-equ}
\left\{
\begin{aligned}
&(\partial_t-\Delta_x )u(t,x,\omega) =0, ~ ~(t,x) \in \mathbb{R}^{1}_+\times \mathbb{R}^n, \quad {\rm a.s.}~\omega,\\
&u(0,x,\omega)=f(x,\omega),~ ~{\rm a.e.\ }\ x \in \mathbb{R}^n,\quad {\rm a.s.}~ \omega,
\end{aligned}
\right.
\end{equation}
with $f \in  L^2(\Omega,  L^2(\mathbb{R}^n)).$ For a.s. $\omega,$ $ f_\omega \in L^2(\mathbb{R}^n)$,~ the convolution
\begin{eqnarray}\label{heat}  u_{f_\omega} (t,x)
=\int_{\mathbb{R}^n} \varphi_t(x-y)f_\omega (y)dy\triangleq(\varphi_t \ast f_\omega) (x)
 \end{eqnarray}
 gives the solution of the classical heat equation with initial condition $f\in L^2(\Omega,  L^2(\mathbb{R}^n)),$
and for each valid $\omega$ the solution $u_{f_\omega} $ belongs to heat-Hardy space $ H^2_{heat}(\mathbb{R}_{+}^{1+n})$ whose space norm is defined as the $L^2(\mathbb{R}^n)$ norm of the non-tangential maximal function $\sup_{(t,y)\in \Gamma^\alpha_x} |u_{f_\omega} (t,y)|,$ where $\Gamma^\alpha_x$ is the orthogonal $\alpha$-cone in $\mathbb{R}_{+}^{1+n}$ with tip at $x\in \mathbb{R}^n.$ The associated stochastic Hardy spaces may be defined.

%\begin{rmk}
%\rm{
% There is a different type of stochastic Dirichlet problem proposed by J.Doob \cite{Doob-sub-harmonic, Bernt-SPDE}.
%}
%\end{rmk}

\begin{Def}\label{particular}(Stochastic harmonic Hardy space)
	The Hilbert space $	L^2(\Omega, h^2(B_1)) $ is a function space that contains all $ h^2(B_1) $-valued random variables $u: B_1 \times \Omega \rightarrow \mathbb{R}$ satisfying
	\begin{equation}
	\begin{aligned}
	&u(\cdot,\omega) \in h^2(B_1), ~{\rm a.s.}~ \omega \in \Omega,\\
	&\mathbb{E}_{\omega}\| u(\cdot,\omega)\|^2_{h^2(B_1)}<\infty.
	\end{aligned}
	\end{equation}
\noindent Analogously, the boundary stochastic Hilbert space $L^2(\Omega, L^2(\partial B_1))$ is defined as:
	\begin{equation}
	\begin{aligned}
L^2(\Omega, L^2(\partial B_1))\triangleq
	\{&f:~\partial B_1\times\Omega\rightarrow\mathbb{R},\\&f(\cdot,\omega) \in L^2(\partial B_1), ~{\rm a.s.}~ \omega \in \Omega,\\
	&\mathbb{E}_{\omega}\| f(\cdot,\omega)\|^2_{L^2(\partial B_1)}<\infty.
	\}.
	\end{aligned}
	\end{equation}
\end{Def}
Stochastic harmonic Hardy space
$L^2(\Omega, h^2(B_1))$ is a Hilbert space with the inner product
\begin{equation}
\langle f,g \rangle_{L^2(\Omega, h^2(B_1))}\triangleq\int_{\Omega}
\int_{\partial B_1}f(t,\omega)g(t,\omega)d\sigma(t)d\mathbb{P}
,~\forall~ f,~g \in L^2(\Omega, h^2(B_1)).
\end{equation}

%It is obvious that $\int_{\partial B_1}f(s,\omega)g(s,\omega)d\sigma(s)$ is well-defined from the definition of $L^2(\Omega, h^2(B_1))$.

Corresponding to the heat equation case a pair of stochastic Hilbert spaces are similarly defined.

\begin{Def}\label{par2}(stochastic heat-Hardy space)
	The Hilbert space $	L^2(\Omega, H^2_{heat}(\mathbb{R}_{+}^{1+n})) $ is a function space that contains all $H^2_{heat}(\mathbb{R}_{+}^{1+n})$-valued random variables $u:  \mathbb{R}_{+}^{1+n}\times \Omega \rightarrow \mathbb{R}$ satisfying
	\begin{equation}
	\begin{aligned}
	&u(\cdot,\omega) \in H^2_{heat}(\mathbb{R}_{+}^{1+n}), ~{\rm a.s.}~ \omega \in \Omega,\\
	&\mathbb{E}_{\omega}\| u(\cdot,\omega)\|^2_{H^2_{heat}(\mathbb{R}_{+}^{1+n})}<\infty.
	\end{aligned}
	\end{equation}
\noindent Analogously, the boundary Hilbert space $L^2(\Omega, L^2(\mathbb{R}^{n}))$ is defined as:
	\begin{equation}
	\begin{aligned}
L^2(\Omega, L^2(\mathbb{R}^{n}))\triangleq
	\{&f:~\mathbb{R}^{n}\times\Omega\rightarrow\mathbb{R},\\&f(\cdot,\omega) \in L^2(\mathbb{R}^{n}), ~{\rm a.s.}~ \omega \in \Omega,\\
	&\mathbb{E}_{\omega}\| f(\cdot,\omega)\|^2_{L^2(\mathbb{R}^{n})}<\infty.
	\}.
	\end{aligned}
	\end{equation}
\end{Def}

As shown in the above contexts the defined stochastic Hardy space are RKHSs. In
spite of the fact that they are isometric, however, the corresponding stochastic $L^2$ spaces are not RKHSs. The mechanism we take for advantage in solving the random boundary value and initial value problems is that in a RKHS, $\mathcal{H}(D),$ the span of the kernels $\{K(\cdot,q)\}|_{q\in \D}$ is dense in $\mathcal{H}(D),$ and the boundary limit functions of the same class $\{K(\cdot,q)\}|_{q\in \D}$ is dense in $L^2(D),$ as well as with equivalent norms. This suggests that we in the $L^2(D)$ do sparse approximation by using the dictionary $\{K(\cdot,q)\}|_{q\in \D}$ and once this is done we $\lq\lq$lift up" the approximating series in the boundary stochastic $L^2$ space to give the stochastic Hardy space approximation in $\mathcal{H}(D).$

We note that the above defined spaces $L^2(\Omega, L^2(\partial B_1))$ and $L^2(\Omega, L^2({\mathbb R}^n))$ in Definitions \ref{particular} and \ref{par2} are boundary stochastic $L^2$-spaces, and $L^2(\Omega, h^2(B_1)) $ and $ L^2(\Omega, H^2_{heat}(\mathbb{R}_{+}^{1+n}))$ are the corresponding stochastic solution spaces. We note that the solution spaces $h^2(B_1)$ and $H^2_{heat}(\mathbb{R}_{+}^{1+n})$ are RKHSs (\cite{Saitoh_RKHS-2016}), briefly denoted as $H^2$ in the sequel.

The idea of POAFD algorithm was initiated in \cite{QSW} and further formulated in \cite{Q2D} and \cite{Qian2018} (also see \cite{CQT}). The work \cite{qu2018} and \cite{qu2019} carry out the POAFD algorithm to the weighted Bergman spaces in the unit disc, and the work \cite{Qu-Sps-represt-dirac} carries out it to the $h^2(B_1)$ and $H^2_{heat}(\mathbb{R}_{+}^{1+n})$ spaces which are all for deterministic signals. The work \cite{QianSAFD} generalizes the AFD and the POAFD methods to random signals. There are two types of generalizations, called SPOAFD1 and SPOAFD2.
Each of them can be applied to the above defined stochastic Hardy spaces.  In this paper, we are based on the SPOAFD methods to solve SPDEs. In the following section we will first outline the deterministic POAFD, and then the SPOAFD algorithms for the self-containing purpose.

\section{ POAFD and SPOAFD Algorithms}\label{Sect-SPOAFD}
For the self-containing purpose we give a revision of POAFD and SPOAFD Algorithms
(see \cite{Qian2018,CQT}, or \cite{qu2018}).

\subsection{POAFD}
In order to perform POAFD the underlying Hilbert space $L^2(\partial \D)$ is assumed to have a subset indexed by elements in $\D,$ denoted $K_q, q\in \D,$ whose span is dense (dense-span property) in $L^2(\partial \D).$ The collection of the reproducing kernels of any RKHS, in particular, has such dense-span property. By this reason, with a little abuse of terminology we call the functions in such dense-span collection by $\lq\lq$kernels".
Besides the density we
further assume the
\emph{boundary vanishing condition} (BVC): For any function $G\in L^2(\partial \D)$ there holds
\begin{eqnarray}\label{BVC}
\lim_{q\to \partial D}|\langle G,E_q\rangle|=0,
\end{eqnarray}
where $E_q$ is the normalized kernel function:
\begin{eqnarray}\label{normalized}
E_q=\frac{K_q}{\|K_q\|}.\end{eqnarray}
 It is necessary to introduce the notion \emph{multiple kernels}. Let $(q_1,\cdots,q_n)$ be an $n$-tuple of parameters in $D.$ Denote by $l(k)$ the multiplicity of $q_k$ in the $k$-tuple $(q_1,\cdots,q_k), 1\leq k\leq n.$
With a little abuse of notation, we define the $k$-th multiple kernel as
\[\tilde{K}_{q_k}=\left[\left(\frac{\partial}{\partial \overline{q}}\right)^{(l(k)-1)}K_q\right]_{q=q_k},\ k=1,2,\cdots,n,\]
where $\frac{\partial}{\partial \overline{q}}$ is interpreted as a directional derivative.
 Multiple kernels are generated from performing
the \emph{pre-orthogonal maximal
 selection principle}: Suppose we already have an $n$-tuple
 $\{q_1,\cdots,q_{n}\},$ possibly with multiplicity. We correspondingly have
 the $n$-tuple of multiple kernels, $\{\tilde{K}_{q_1},\cdots,\tilde{K}_{q_{n}}\}.$
 Applying the G-S orthonormalization process consecutively, we obtain
 an equivalent $n$-orthonormal system,
 $\{E_1,\cdots,E_{n}\},$ where $E_1=E_{q_1}.$
 For any given $G$ in the Hilbert space we now investigate whether there exists a $q_{n+1}\in D$ that gives rise to the supreme value
 \[ |\langle G,E_{n+1}\rangle|=\sup\{|\langle G,E_{n+1}^q\rangle| \ :\ q\in {D},~ q\ne q_1,\cdots,q_{n}\},\]
where the finiteness of the supreme is guaranteed by the Cauchy-Schwarz inequality, and
$E_{n+1}^q$ be such that $\{E_1,\cdots,E_{n},E_{n+1}^q\}$ is the G-S orthonormalization
of $\{\tilde{K}_{q_1},\cdots,\tilde{K}_{q_{n}}, K_q\},$ and $E_{n+1}$ be such that $\{E_1,\cdots,E_{n},E_{n+1}\}$ is the G-S orthonormalization
of $\{\tilde{K}_{q_1},\cdots,\tilde{K}_{q_{n}}, \tilde{K}_{q_{n+1}}\}.$ In fact, since $q$ is distinct from the proceeding $q_1,\cdots,q_{n},$ $E_{n+1}^q$ is given by
 \begin{eqnarray}\label{GS}
 E_{n+1}^q
 =\frac{K_q-\sum_{k=1}^{n}\langle K_q,E_k\rangle E_k}
 {\|K_q-\sum_{k=1}^{n}\langle K_q,E_k\rangle E_k\|}.\end{eqnarray}
Under BVC a compact argument concludes that
 there exists a point $q_{n+1}\in D$ and $q^{(l)}, l=1,2,\cdots,$
  such that $q^{(l)}$ are all different from $q_1,\cdots,q_{n},$
  $\lim_{l\to \infty}q^{(l)}=q_{n+1},$ and
 \begin{eqnarray}\label{max} \lim_{l\to\infty}|\langle G,E_{n+1}^{q^{(l)}}\rangle|=
\sup\{|\langle G,E_{n+1}^q\rangle| \ :\ q\in {D}, q\ne q_1,\cdots,q_{n}\}=|\langle G,E_{n+1}\rangle|,\end{eqnarray}
 where, by a Taylor expansion argument involving the Lagrange remainder, we have
 \begin{eqnarray}\label{GSmultiple}
 E_{n+1}
 =\frac{\tilde{K}_{q_{n+1}}-\sum_{k=1}^{n}\langle \tilde{K}_{q_{n+1}},E_k\rangle E_k}
 {\sqrt{\|\tilde{K}_{q_n}\|^2-\sum_{k=1}^{n}|\langle \tilde{K}_{q_{n+1}},E_k\rangle|^2}}\end{eqnarray}
 (see \cite{Qian2018,CQT}, or \cite{qu2018}).
Iteratively applying the above process to $G,$ we obtain
 a sequence of maximally selected $\{q_k\}_{k=1}^\infty,$ and  has
\begin{eqnarray}\label{converge} G=\sum_{k=1}^\infty \langle G,E_k\rangle E_k.\end{eqnarray}
We will call $\{E_k\}$ the consecutive Gram-Schmidt orthogonalization of the multiple kernels corresponding to maximally selected $\{q_k\},$ the latter being possibly with multiples. It is noted that the outcome system $\{E_k\}$ is not necessarily a basis, it, however, is made for, and can well express the given signal. For the POAFD expansion we have in the space norm sense the same convergent rate $O(n^{-1/2})$ as that for the K-L expansion (\cite{Q2D,Qian2018,CQT,LPS}).

Next, we incorporate probability and formulate stochastic POAFD 1 and 2 (SPOAFD1 and 2). For the stochastic Dirichlet problem the set $D$ stands for the unit ball $B_1\in {\mathbb R}^{1+n}$ or the upper half space $\mathbb{R}_{+}^{1+n},$ while for the stochastic Cauchy initial value case $D$ stands for the half space $\mathbb{R}_{+}^{1+n}.$  By involving such domain $D$ we adopt the setting $\N=L^2(\Omega,  L^2(\partial D))$  with any underlying probability space $(\Omega, \mathcal{F},\mathbb{P}).$ We will be related to the stochastic Hardy spaces $L^2(\Omega, H^2),$ where $H^2$ stands for either  $h^2(B_1),$ or $h^2(\mathbb{R}_{+}^{1+n}),$ and $H^2_{heat}(\mathbb{R}_{+}^{1+n}),$ depending on the individual problem. Their reproducing kernels are uniformly denoted as $K(\cdot,q)=K_q(\cdot),$ whose span is dense in the concerned $L^2(\Omega, H^2),$ as well in $L^2(\Omega, L^2(\D)).$ The stochastic Hardy spaces are respectively imbedded,
through the non-tangential boundary limits, into the $\N$ spaces, namely $\N=L^2(\Omega, L^2(\partial B_1))$ or $\N=L^2(\Omega, L^2({\mathbb R}^n)).$

\subsection{SPOAFD1}

The SPOAFD1 algorithm is designed for a stochastic signal $f(t,\omega)\triangleq f_\omega (t),$ being expressible as $ f_\omega (t)\triangleq\tilde{f}(t)+r_\omega(t),$ where $\tilde{f}(t)\triangleq\mathbb{E}(f_\omega(t))$ and
$r_\omega(t)\triangleq f_\omega(t)-\tilde{f}(t),$ and the error random signal $r_\omega(t)$ is assumed to be of small $\N$-norm. Obviously, $\mathbb{E}(r_\omega(t))=0$ for a.e. $t\in D$. SPOAFD1 is designed to analyze a given stochastic signal $f$ by means of the POAFD expanding orthonormal system $\{E_k\}_k$ of the deterministic signal $\tilde{f}.$

The SPOAFD1 method to solve the stochastic boundary value or initial value problems consists of three steps. Below as examples we work on the stochastic Dirichlet problem on the unit ball and on the stochastic Cauchy problem on the upper-half space.
  The first step is to compute out $\tilde{f}.$  The second step is to apply POAFD to the deterministic signal $\tilde{f}$ with the expansion
\begin{eqnarray}\label{CS1}
\tilde{f}=\sum_{k=1}^\infty \langle \tilde{f},E_k\rangle E_k
\end{eqnarray}
valid on the boundary, where $\{E_k\}$ is consecutively the multiple G-S orthogonalizations of the multiple kernels corresponding to $\{q_k\}.$  The third step is to term by term solve the problems and add them together by
invoking boundedness of the Poisson kernel convolution operators. The semigroup property of the spherical Poisson kernel then can $\lq\lq$lift up" the expansion on the sphere and, owing to the isometry between the $h^2(\D)$ and the $L^2(\partial \D)$ spaces, give
 (see \S 4 of \cite{Qu-Sps-represt-dirac})
\[ u_{\tilde{f}}(x)=\sum_{k=1}^\infty \langle \tilde{f},E_k\rangle u_{E_k}(x)\triangleq\sum_{k=1}^\infty \langle \tilde{f},E_k\rangle P(E_k)(x)=\sum_{k=1}^\infty c_kP_{|x|q_k}(x'),\qquad x\in B_1,\]
where $P(E_k)$ is the classical solution of the Dirichelet problem with boundary data $E_k,$ $ x=|x|x^\prime,~q_k=|q_k|q_k^\prime$ and $\{c_k\}_k$ are certain constants
(see the later half of the proof of Theorem \ref{5} below).
Similarly, for the stochastic Cauchy initial value problem, the SPOAFD1 expansion together with the lifting up based on the semigroup property gives rise to the solution
(see \S 3.2 of \cite{Qu-Sps-represt-dirac})
\[ u_{\tilde{f}}(x)=\sum_{k=1}^\infty \langle \tilde{f},E_k\rangle u_{E_k}(x)\triangleq\sum_{k=1}^\infty \langle \tilde{f},E_k\rangle P_t\ast E_k(\underline{x})=\sum_{k=1}^\infty d_kh_{t+q_k}(\underline{x}),\qquad x\in \mathbb R^{n+1}_+,\]
 where $x=t+\underline{x}, q_k=s_k+\underline{q_k},$ and $\{c_k\}_k$ are certain constants. Note that the convergence rates of the solution series are as the same as those expanding the stochastic boundary or the initial values.
Below, like $u_{\tilde{f}},$ we will use the self-explanatory notations $u_f$ and $u_r.$
Denote
 \begin{eqnarray}\label{d} d_{u_{f_\omega}}(x)=u_{f_\omega}(x)-\sum_{k=1}^\infty \langle {f_\omega},E_k\rangle u_{E_k}(x).\end{eqnarray}
 Since $\|r\|^2_\N=\|{\rm var}f\|_{L^1(\partial \D)},$ the results
 of \cite{QianSAFD} imply

\begin{Thm}
\[ \|u_{f_\omega}-u_{\tilde{f}}\|_{\N}\leq \|{\rm var}f\|_{L^1(D)},\]
and
\[ \|d_u\|_{\N}\leq \|{\rm var}f\|_{L^1(\partial \D)}-\sum_{j=1}^{\infty}\mathbb{E}|\langle r_{\omega},E_j\rangle|^2,
\]
and, in particular,
\[ \|d_u\|_{\N}\leq \|{\rm var}f\|_{L^1(\partial \D)}.\]
\end{Thm}

 The theorem concludes that SPOAFD1 solves the stochastic boundary and initial problems up to an error dominated by the $L^1(\partial \D)$-average of the variation of $f.$

\subsection{SPOAFD2}
SPOAFD2, or, abbreviated as SPOAFD, is the main method introduced by this paper, because through it we gain  a.s. convergence of the involved series expansions in the $H^2(\D)$ and the $L^2(\partial \D)$ spaces.  SPOAFD is built upon the following \emph{statistical maximal selection principle} (\emph{SMSP}).

\begin{Thm}
	Assume that $f\in L^2(\Omega, L^2(\partial \D)),$ where with respect to the dictionary $\{K_q\}$ the space $L^2(\partial \D)$ satisfies BVC. There exists $q_1 \in D$ such that
	\begin{equation}
	q_1=\arg\max_{q \in D}\mathbb{E}|\langle f_\omega,E_q\rangle|^2 ,
	\end{equation}
	where $E_q$ is the normalized reproducing kernel of the Hardy space $H^2$ defined as in (\ref{normalized}) for either the Dirichlet boundary value or the Cauchy initial value problem.
\end{Thm}

\noindent\textbf{Proof}:
 By invoking the main result of \cite{QianSAFD} it is sufficient to prove the \emph{statistical boundary vanishing condition} property, i.e.\\
\begin{equation}
\lim_{q\rightarrow \partial D}\mathbb{E}|\langle f_\omega,E_q\rangle|^2 =0.
\end{equation}
Since  \[  f_\omega \in L^2(\partial \D),~{\rm a.s.}~ \omega \in \Omega, \]
we have
\begin{equation}\label{h2-rpdc-ker-BVC}
\lim_{q\rightarrow \partial D}|\langle f_\omega,E_q\rangle|^2 =0,~{\rm a.s.}~ \omega \in  \Omega,
\end{equation}
as proved in \cite{Qu-Sps-represt-dirac}.
Noting that, for $q\to\partial \D$ the functions $|\langle f_\omega,E_q\rangle|^2$ in $\omega$ is dominated uniformly in $q$ by an integrable function in $L^1(\Omega,d\mathbb{P}):$
\begin{equation}
\begin{aligned}
|\langle f_\omega,E_q\rangle|^2& \leq \|f_\omega\|^2_{H^2}\|E_q\|^2_{H^2}\\&=\|f_\omega\|^2_{H^2} \in L^1(\Omega,d\mathbb{P}),
\end{aligned}
\end{equation}
by \textit{Lebesgue's dominated convergence theorem},~we have
\begin{equation}
\lim_{|q|\rightarrow \partial D}\mathbb{E}|\langle f_\omega,E_q\rangle|^2 =0.
\end{equation}
By employing the usual compact argument we obtain a maximal selection of $q\in D.$ The proof is complete.$\eop$

 The availability of SMSP implies that SPOAFD2 can be performed in $ L^2(\Omega, L^2(B_1)) $ and $L^2(\Omega,  L^2(\mathbb{R}^n)).$ That is,

\begin{Thm}\label{see} (\cite{QianSAFD})
	Suppose $ f \in \N= L^2(\Omega, L^2(\D)).$  Assume that $ \{q_j\}_j $
	are selected under SMSP,
	\begin{equation}\label{Statistical-average}
	q_l=\arg\max_{q \in D}\mathbb{E}|\langle f_\omega,E_l\rangle|^2 ,
		\end{equation}
	where $\{E_k\}$ is the consecutive G-S orthogonalization of the multiple kernels corresponding to the maximally selected $\{q_k\}.$  Then there holds
	\begin{equation}\label{SPOAFD-h1}
	{f_\omega}\stackrel{\mathcal{N}}{=}\sum_{k=1}^{\infty}\langle f_\omega,E_k\rangle {E_k}.
	\end{equation}
\end{Thm}

As a consequence of the above theorem, we have
\begin{Thm}\label{5}
Under the SPOAFD2 expansion (\ref{SPOAFD-h1}) there holds
\begin{equation}\label{SPOAFD-h2}
	u_{f_\omega}\stackrel{\mathcal{N}}{=}\sum_{k=1}^{\infty}\langle f_\omega,E_k\rangle u_{E_k}=\sum_{k=1}^{\infty}c_k(\omega) \tilde{K}_{q_k},
	\end{equation}
where $\{\tilde{K}_{q_k}\}$ are the multiple kernels whose consecutive G-S orthogonalizations are $\{E_k\}$ computed as in (\ref{GSmultiple}), $\{q_k\}$ are the parameter sequence consecutively obtained through applying SMSP to $f_\omega,$  $\{c_k(\omega)\}$ are accordingly asserted coefficients depending on $\omega.$
\end{Thm}

\noindent {\bf Proof}:~The POAFD convergence result amounts, for each $\omega,$ possibly excluding an event of probability zero,
\[ f_\omega \stackrel{L^2(\partial D)}{=} \sum_{k=1}^\infty \langle f_\omega, E_k \rangle E_k.\]
Since the lifting up  operators $f_\omega \mapsto u_{f_\omega}$ for solving the classical problems is either isometric or bounded (for the heat kernel case), we have, for a.s. $\omega\in \Omega,$
\[ u_{f_\omega} \stackrel{H^2(D)}{=} \sum_{k=1}^\infty \langle f_\omega, E_k \rangle u_{E_k}.\]
Since the Hardy space projection $u_f=P_{H^2}(f): L^2(\Omega,L^2(\partial D))\to L^2(\Omega,H^2(D))$ is linear and isometric, we have, for any fixed $n,$
\begin{eqnarray*}
 \|u_{f_\omega}-\sum_{k=1}^n \langle f_\omega, E_k \rangle u_{E_k}\|_{L^2(\Omega,H^2(D))}&\leq &
  \|f_\omega-\sum_{k=1}^n \langle f_\omega, E_k \rangle E_k\|_{L^2(\Omega,L^2(\partial D))}. \end{eqnarray*}
 By invoking the SPOAFD2 convergence result in \cite{QianSAFD} the last quantity tends to zero as $n\to \infty.$ Therefore, (\ref{SPOAFD-h2}) holds.
 Next, we deduce the coefficients $c_k(\omega).$

 Denote, for every positive integer $n$,~$\mathcal{E}_n=(E_1,\cdots,E_n),~ \mathcal{K}_n=(\tilde{K}_{q_1},\cdots,~\tilde{K}_{q_n}),~ \mathcal{C}_n=(c_1(\omega),\cdots,c_n(\omega)),$ and $ \mathcal{F}_n=(\langle f_\omega ,E_1\rangle, \cdots, \langle f_\omega ,E_n\rangle ).$
 Then there follows
 \[
 {\mathcal{K}_n}^t=\mathcal{A}_n{\mathcal{E}_n}^t,\]
 where the superscript $\lq\lq$t" stands for transpose of the matrix, and
 $\mathcal{A}_n=(a_{ij})_{1\leq i,j\leq n}$ is a lower triangle $n\times n$ matrix, where
 \[ a_{ij}=
 \left\{\begin{array}{ll}
 E_{q_i}(q_j)=
 \frac{K(q_j,q_i)}{\|K_{q_i}\|},~i\ge j,\\
 0,~i<j.
 \end{array}\right.\]
  Hence,
 \[ {\mathcal{E}_n}^t=\mathcal{A}_n^{-1}{\mathcal{K}_n}^t,\quad {\rm and}\quad
 \mathcal{C}_n=\mathcal{F}_n \mathcal{A}_n^{-1}.\]
The proof is complete.$\eop$

Theorem \ref{5} shows that with $f$ being in the Bochner space $L^2(\Omega,L^2(\partial D))$ the stochastic Dirichlet boundary value problem and the stochastic Cauchy initial value problem are solvable in their corresponding Bochner type stochastic Hardy spaces $L^2(\Omega,H^2)$ with explicit sparse representations.

We note that theoretically the introduced method has the standard convergence rate $\frac{M}{\sqrt{n}}$ (see \cite{Qian2018,CQT} or \cite{KLvsAFD}).

\section{Examples of Random Boundary Value Problem}
\label{Sect-experiment}
In this section we first present the algorithm of SPOAFD2 (SPOAFD in short), which is applied to numerical simulations of two types of random PDEs described by the Laplace equation and the heat equation. Denote by $\{B_1,\cdots,B_j\}$ the G-S orthonormalization of the maximally selected parameterized multiple kernels under the SPOAFD program.\\

\begin{tabular}{l}
	\hline
	{\bfseries Algorithm.} SPOAFD \\
	\hline
	\;{\bfseries Input:}~~$\quad$original function $f(t,\omega)$\\
	\;{\bfseries Output:} ~solution $u_n$ of random equation \\
	\;1:$\qquad\quad$  ~{ Initialize} $u_n=0,~j=0,$ error=10, ${\rm energy}_\omega=0.$\\
	\;2:$\qquad\quad$ ~~{\bfseries While} error $\geq$ $10^{-4}$ {\bfseries do}\\
	\;3:$\qquad\quad\qquad$ ~~$j\leftarrow j+1$\\
	\;4:$\qquad\quad\qquad$ ~~$a_j\leftarrow$ $\arg\max\limits_{b\in D} \mathbb{E}_\omega|\langle f_\omega,B_j^b\rangle|^2$\\
	\;5:$\qquad\quad\qquad$ ~~$u_n\leftarrow u_n+\langle f_\omega,B_j^{a_j}\rangle B_j^{a_j}$\\
	\;6:$\qquad\quad\qquad$ ~~${\rm energy}_\omega \leftarrow {\rm energy}_\omega+|\langle f_\omega,B_j^{a_j}\rangle|^2$\\
	\;7:$\qquad\quad\qquad$ ~~error~$\leftarrow\frac{\|f\|_{\mathcal{N}}^2-\mathbb{E}({\rm energy}_\omega)}{\|f\|_{\mathcal{N}}^2}$\\
	\;8:$\qquad\quad$ ~~{\bfseries End while}\\
	\hline\par
\end{tabular}

\bigskip
As already noted that the expansions using the parameterised fundamental solutions have the convenience of lifting up to get the virtual solutions of the problems. For such purpose below we only expand the random boundary data by the respective dictionaries containing the needed fundamental solutions. Without such convenience random boundary data can also be expanded by other basis or dictionary elements, such as the Fourier-Hermite polynomials and the eigenfunctions of K-L decomposition. A separate paper, \cite{KLvsAFD}, will devote to comparison between the existing stochastic expansions and SPOAFD, without involving applications to SPDEs.

%\begin{tabular}{l}
%	\hline
%	{\bfseries Algorithm 1.} SPOAFD1 \\
%	\hline
%	\;{\bfseries Input:}~~$\quad$original function $f(t,\omega)$\\
%	\;{\bfseries Output:} ~solution $u_n$ of the random equation \\
%	\;1:$\qquad\quad$  ~{ Initialize} $u_n=0,~j=0,$ error=10, energy=0.\\
%	\;2:$\qquad\quad$  ~~Computing Mathematical expectation $\tilde{f}(t):=\mathbb{E}(f(t,\cdot)),~\forall~t.$ \\
%	\;3:$\qquad\quad$ ~~{\bfseries While} error $\geq$ $10^{-4}$ {\bfseries do}\\
%	\;4:$\qquad\quad\qquad$ ~~$j\leftarrow j+1$\\
%	\;5:$\qquad\quad\qquad$ ~~$a_j\leftarrow$ arg $\max\limits_{b\in D} |\langle \tilde{f},B_j^b\rangle|$\\
%	\;6:$\qquad\quad\qquad$ ~~$u_n\leftarrow u_n+\langle \tilde{f},B_j^{a_j}\rangle B_j^{a_j}$\\
%	\;7:$\qquad\quad\qquad$ ~~energy~$\leftarrow$~energy+$|\langle \tilde{f},B_j^{a_j}\rangle|^2$\\
%	\;8:$\qquad\quad\qquad$ ~~error~$:=\frac{\|\tilde{f}\|^2-{\rm energy}}{\|\tilde{f}\|^2}$\\
%	\;9:$\qquad\quad$ ~~{\bfseries End while}\\
%	\hline\\
%\end{tabular}

\begin{figure}[H]
	\begin{minipage}[c]{0.3\textwidth}
		\centering
		\includegraphics[height=4.5cm,width=5cm]{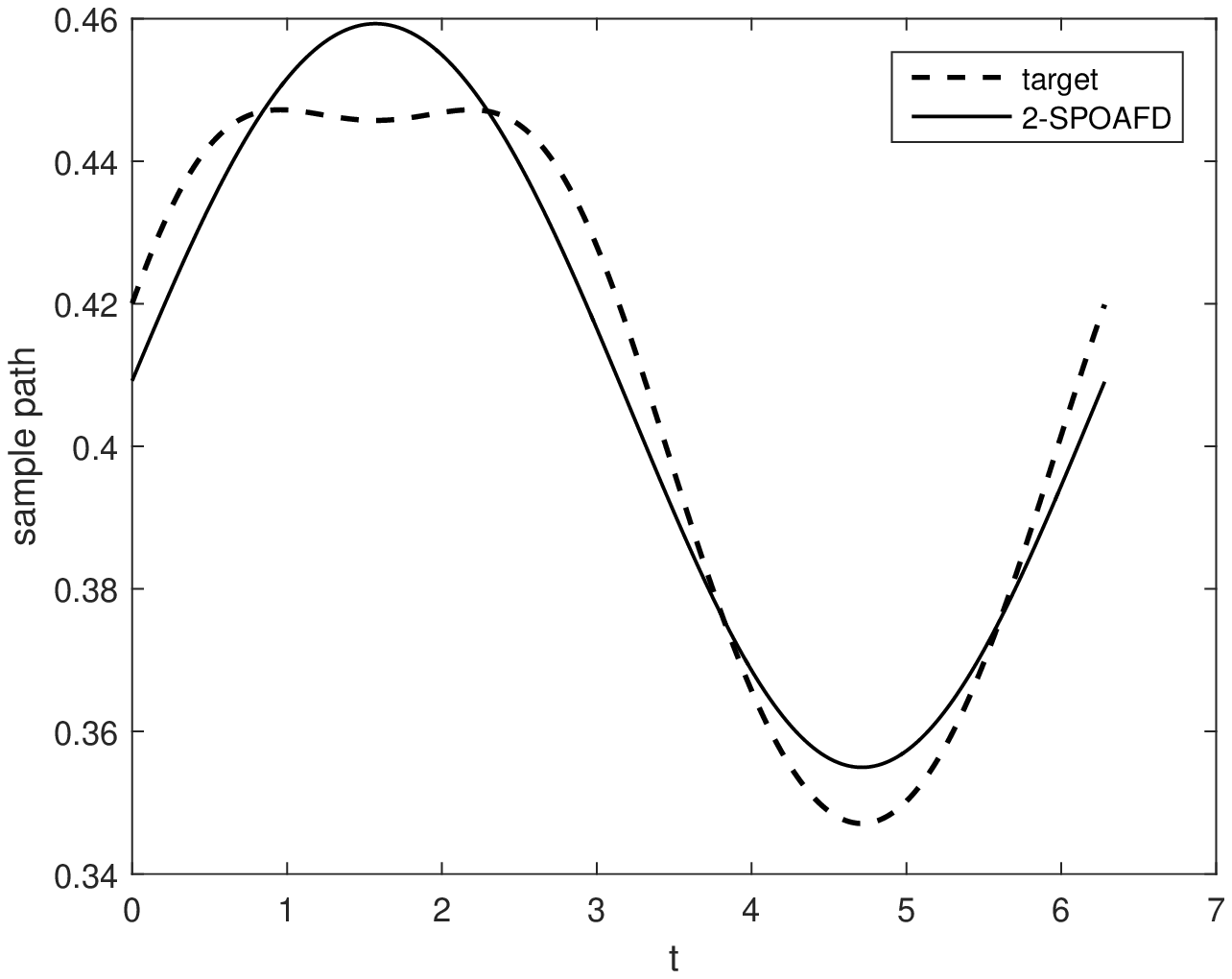}
		\centerline{{\scriptsize  2 partial sum}}
%\centerline{{\scriptsize  relative error: 3.6902$\times 10^{-4}$}}
	\end{minipage}
	\begin{minipage}[c]{0.3\textwidth}
		\centering
		\includegraphics[height=4.5cm,width=5cm]{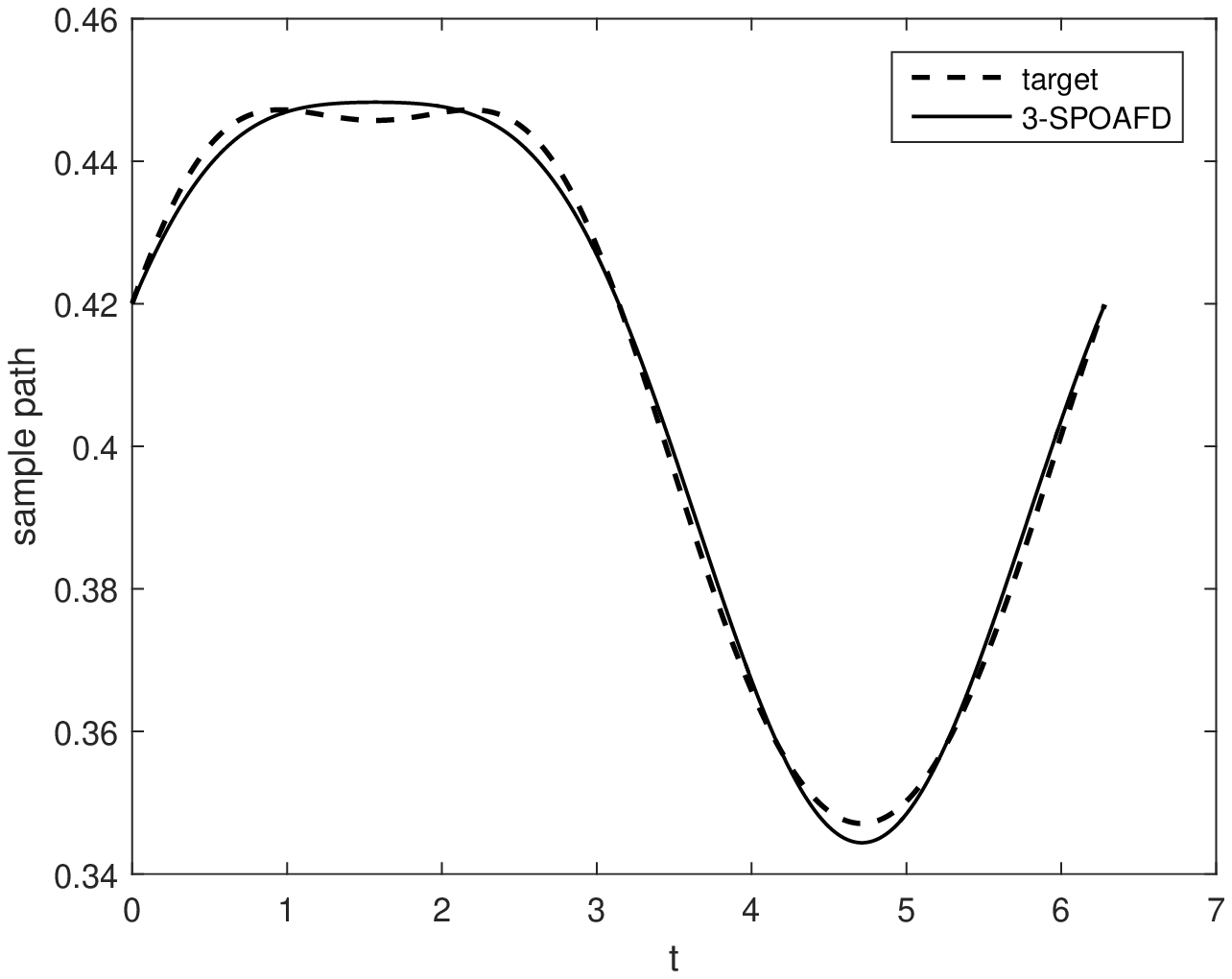}
		\centerline{{\scriptsize 3 partial sum}}
%\centerline{{\scriptsize  relative error: 2.2177$\times 10^{-5}$}}
	\end{minipage}
	\begin{minipage}[c]{0.3\textwidth}
		\centering
		\includegraphics[height=4.5cm,width=5cm]{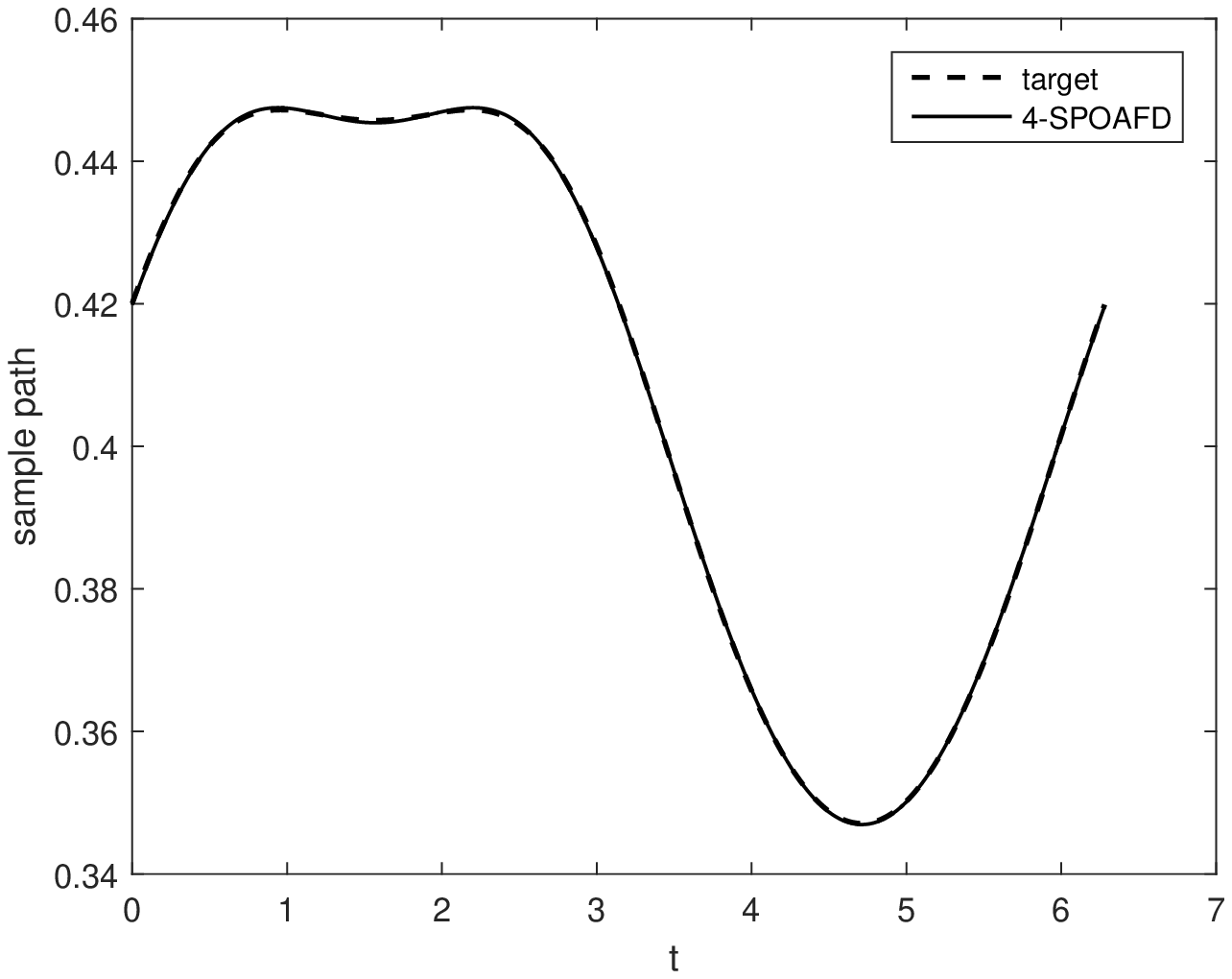}
		\centerline{{\scriptsize  4 partial sum}}
%\centerline{{\scriptsize  relative error: 1.8945$\times 10^{-7}$}}
	\end{minipage}
	
	\caption{$X=0$}
\label{figure1}
\end{figure}

\begin{figure}[H]
	\begin{minipage}[c]{0.3\textwidth}
		\centering
		\includegraphics[height=4.5cm,width=5cm]{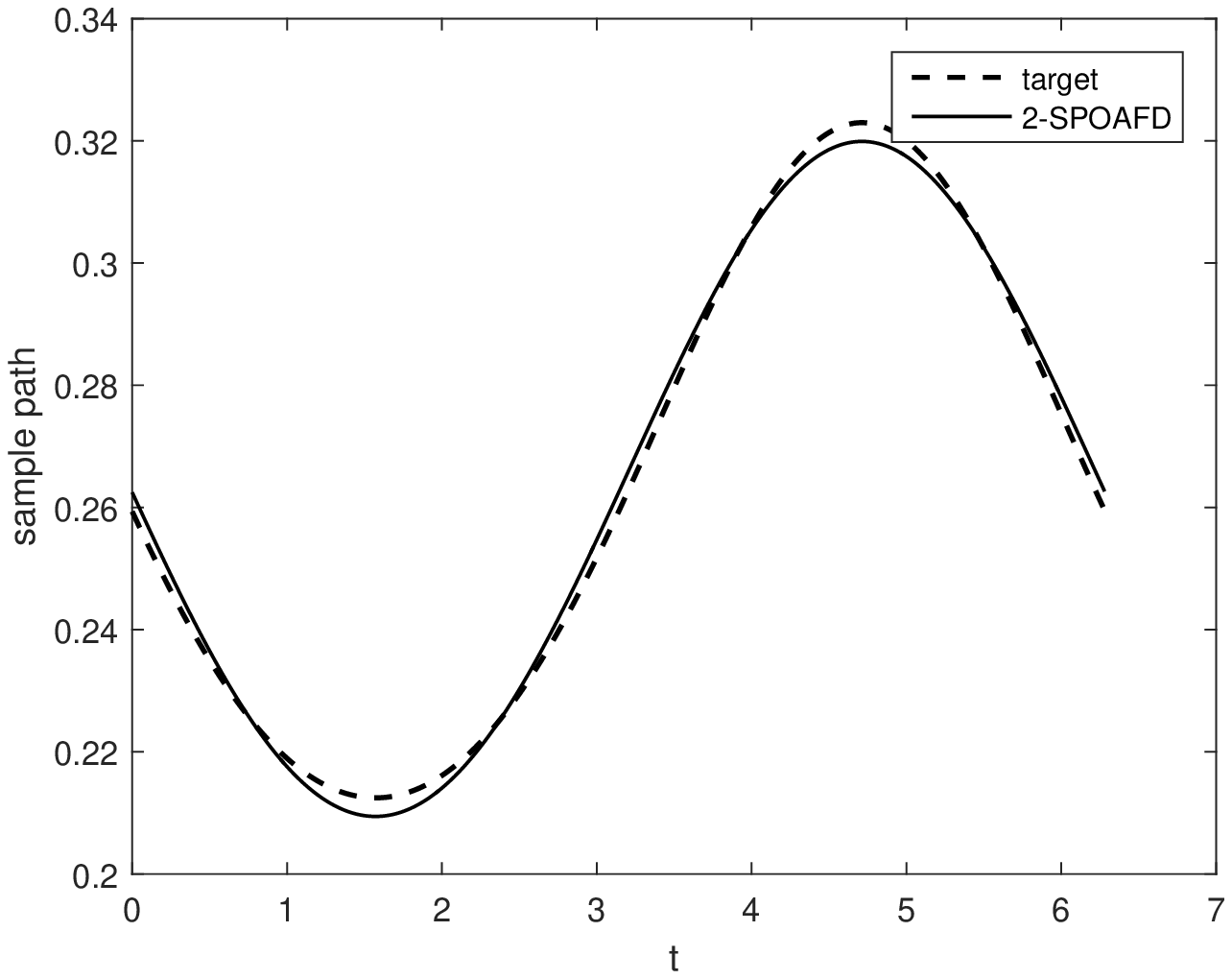}
		\centerline{{\scriptsize  2 partial sum}}
	\end{minipage}
	\begin{minipage}[c]{0.3\textwidth}
		\centering
		\includegraphics[height=4.5cm,width=5cm]{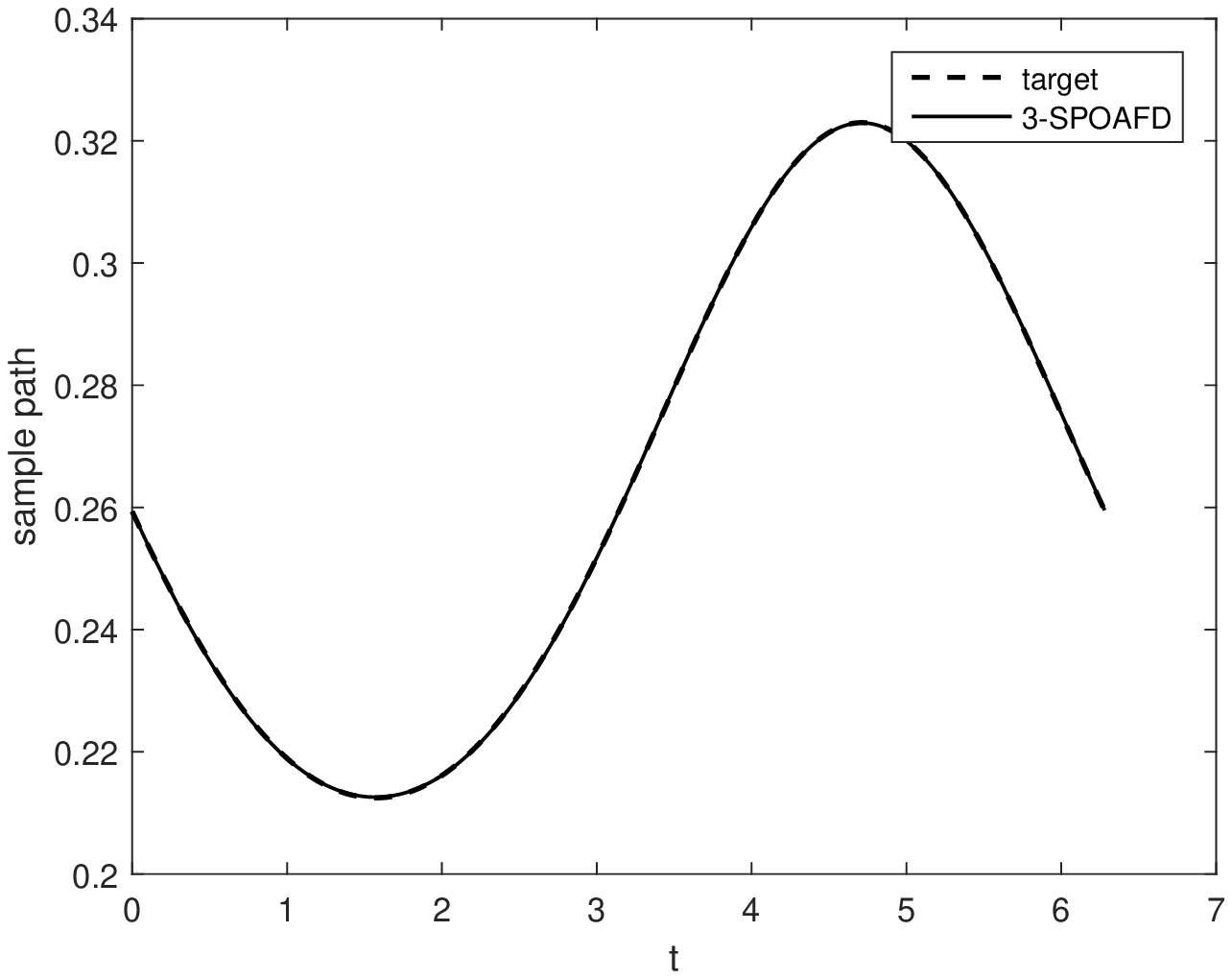}
		\centerline{{\scriptsize  3 partial sum}}
	\end{minipage}
	\begin{minipage}[c]{0.3\textwidth}
		\centering
		\includegraphics[height=4.5cm,width=5cm]{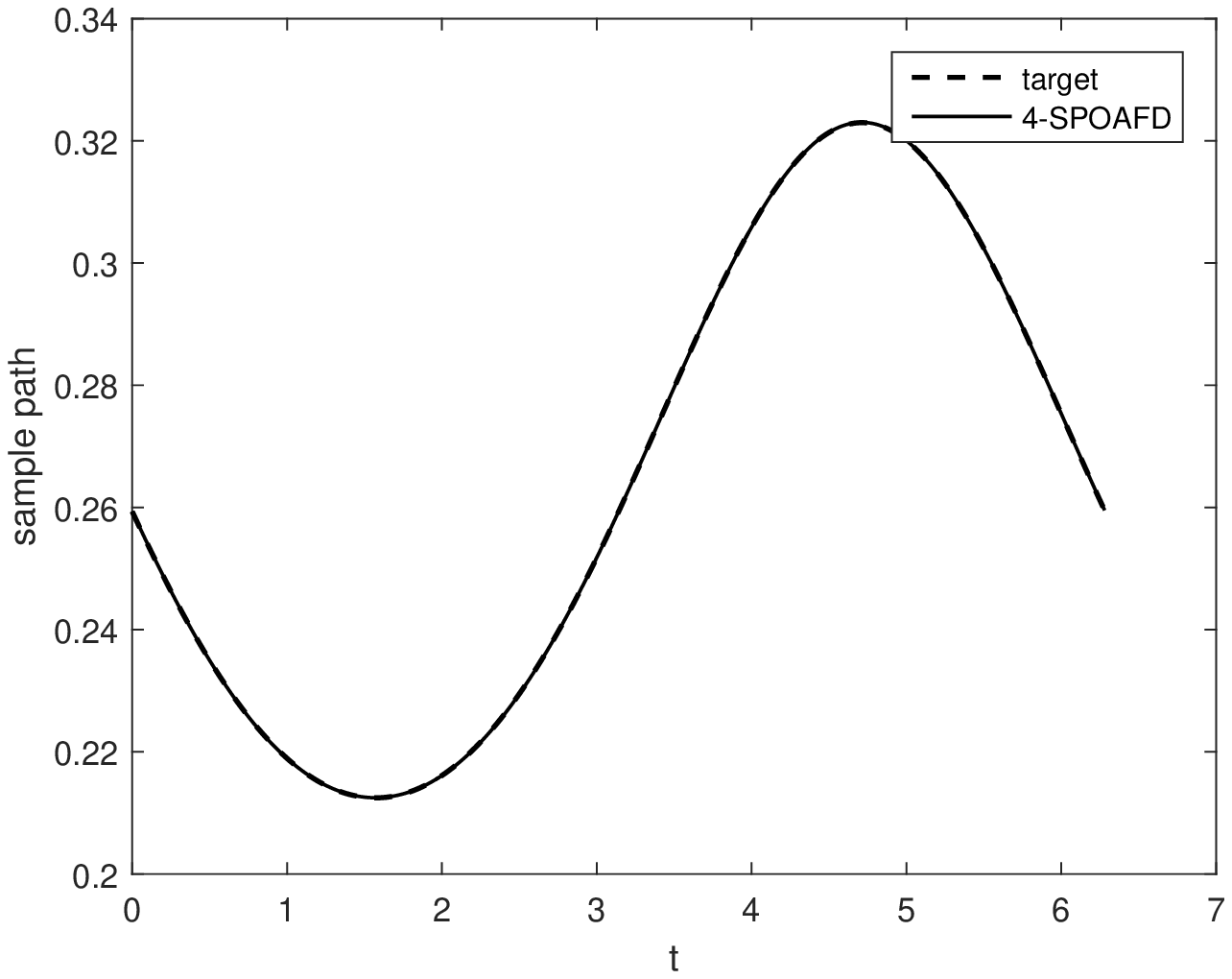}
		\centerline{{\scriptsize  4 partial sum}}
	\end{minipage}

	\caption{$X=-\pi$}
\label{figure2}
\end{figure}

\begin{figure}[H]
	\begin{minipage}[c]{0.3\textwidth}
		\centering
		\includegraphics[height=4.5cm,width=5cm]{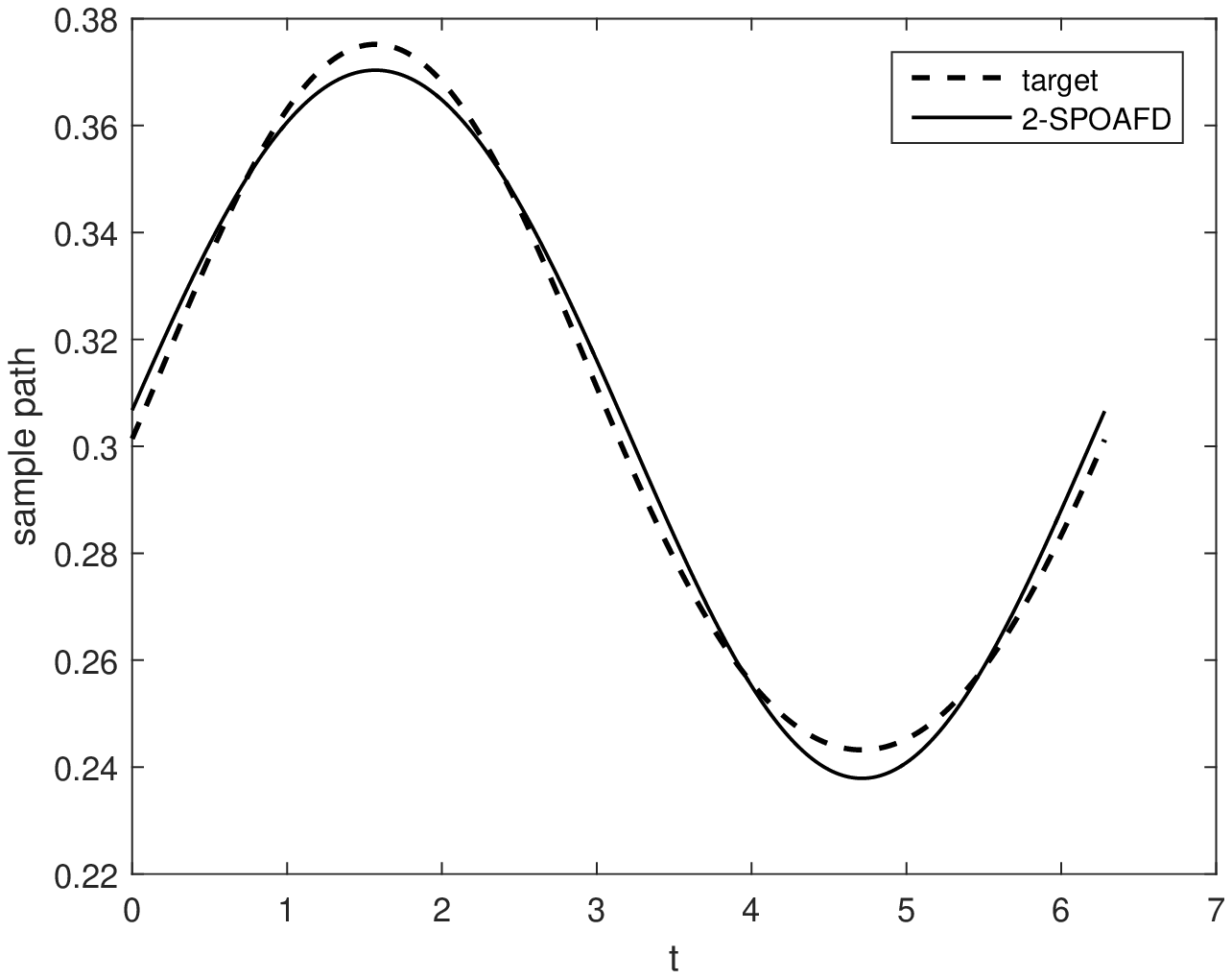}
		\centerline{{\scriptsize 2 partial sum}}
	\end{minipage}
	\begin{minipage}[c]{0.3\textwidth}
		\centering
		\includegraphics[height=4.5cm,width=5cm]{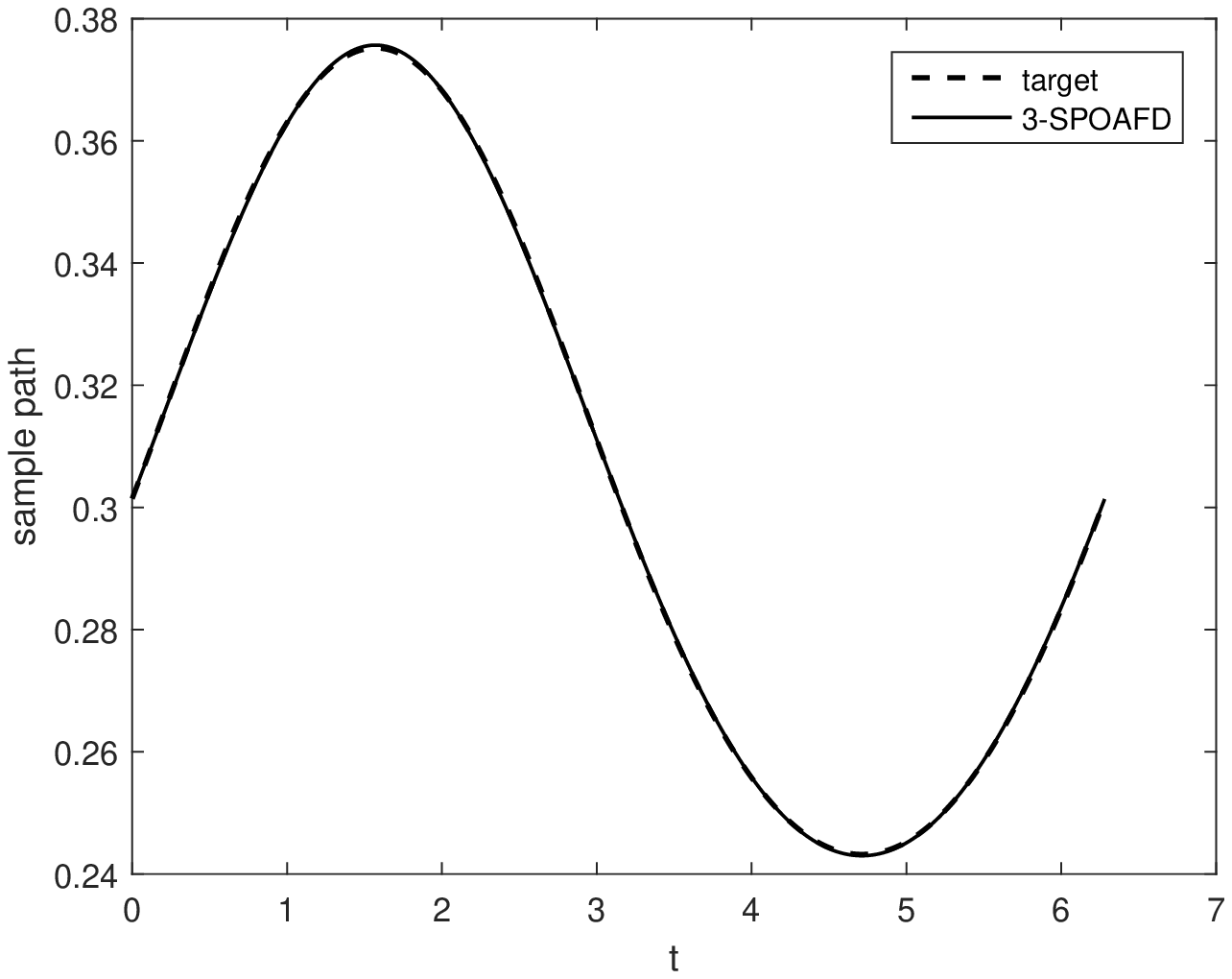}
		\centerline{{\scriptsize 3 partial sum}}
	\end{minipage}
	\begin{minipage}[c]{0.3\textwidth}
		\centering
		\includegraphics[height=4.5cm,width=5cm]{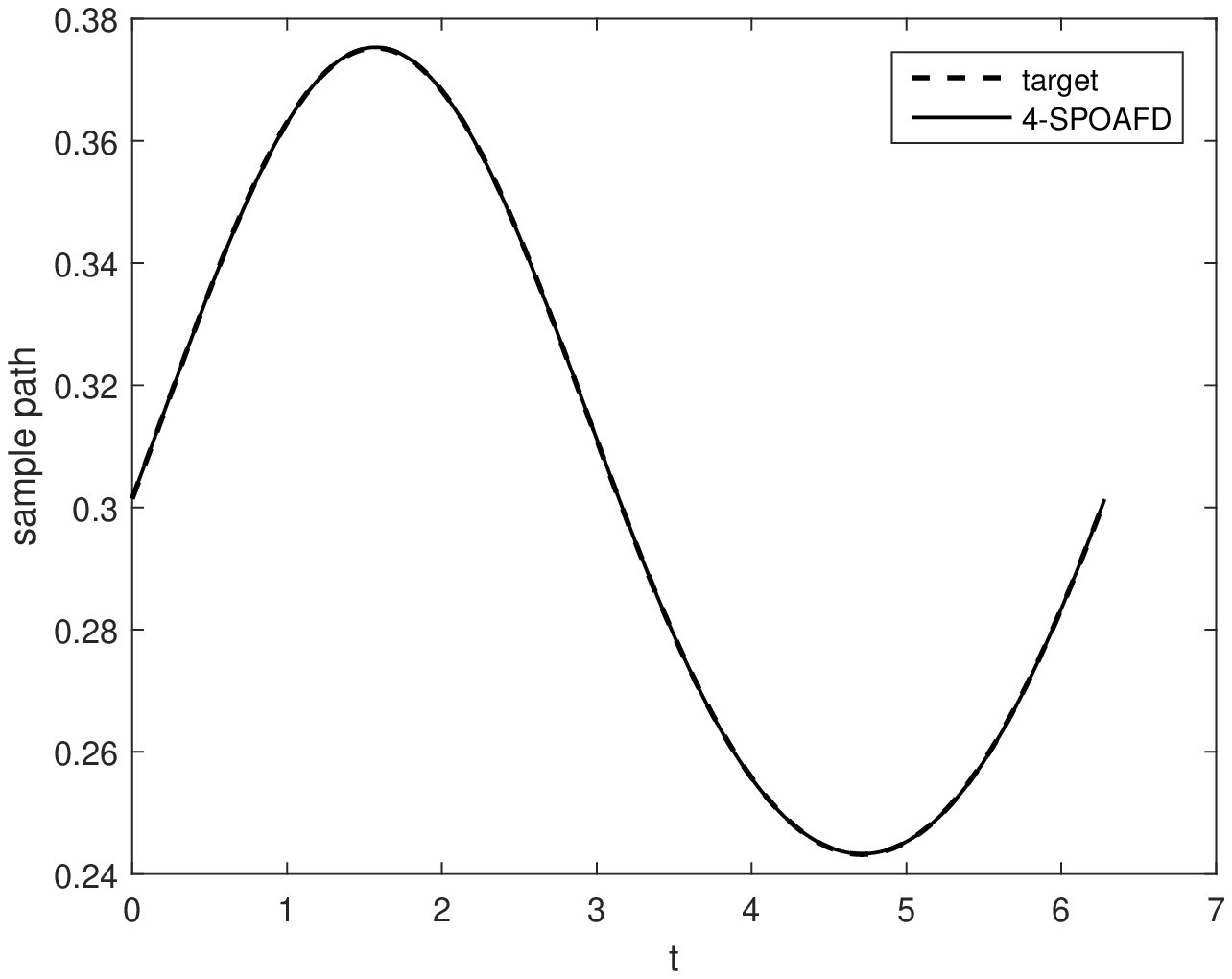}
		\centerline{{\scriptsize 4 partial sum}}
	\end{minipage}

	\caption{$X=2.4504$}
\label{figure3}
\end{figure}

\begin{Exmp}\label{Exmp-Lplc}

Consider the Laplace equation in the unit disk $\mathbb{D}\triangleq\{z\in \mathbb{C}:~|z|<1 \}$ with the random Dirichlet boundary condition:
\begin{equation}
\left\{
\begin{aligned}
&\Delta u(z,\omega)=0,~ \forall~z \in\mathbb D,~{\rm a.s.} ~\omega \in \Omega,\\
 &u(e^{it},\omega)=1/\sqrt{5+(\sin t-X)^2},~ ~ e^{it} \in \partial \mathbb{D},~{\rm a.s.},
\end{aligned}
\right.
\end{equation}
where $t\in [0,2\pi]$ and $X=X(\omega)$ is a random variable with the density function
$p(s)=\frac{1}{m}\alpha(s),$ where $$\alpha(s)=\begin{cases} e^{-\frac{1}{(s-\pi)^2+1}}, & 0\leq s \leq\pi, \\ e^{-\frac{1}{(s+\pi)^2+1}}, & -\pi\leq s<0, \end{cases}$$
and $m=\int_{-\pi}^{\pi}\alpha(s)ds.$
\end{Exmp}
\noindent
%We note that
%\[ \mathbb{E}(u(e^{it},\cdot))=e^{\cos(t-\pi)}=:\tilde{f}(e^{it}).\]
   Since the boundary data is a bi-variate function of $t$ and $X,$ we can make use the distribution of the random variable $X.$ We have three groups of figures respectively for visual effects of the experiments at those $\omega$ such that $X(\omega) =0, -\pi,$ and $2.4504,$ which can be shown in Figure \ref{figure1}, \ref{figure2} and \ref{figure3}. The relative errors in Table \ref{Lplc-table1},~\ref{Lplc-table2},~\ref{Lplc-table3} are computed according to the formulas given in the relevant algorithms.
%The parameters selected according to the respective maximal selection principle are

%1

\begin{table}[H]
\scriptsize
\centering
	\begin{tabular}{c|c|c|c}
		\toprule
		SPOAFD&3.6902$\times 10^{-4}$ & 2.2177$\times 10^{-5}$ & 1.8945$\times 10^{-7}$ \\
%		\midrule
%		numerical KL&4.7292$\times 10^{-5}$& 6.0649$\times 10^{-6}$& 5.5206$\times 10^{-8}$ \\
		\bottomrule
	\end{tabular}%
	\centering
	\caption{\scriptsize relative error}
	\label{Lplc-table1}%
\end{table}%

%2

\begin{table}[H]
\scriptsize
	\begin{tabular}{l|l|l|l}
		\toprule
		SPOAFD&6.7418$\times 10^{-5}$& 6.2401$\times 10^{-8}$& $6.7669\times 10^{-9}$\\
%		\midrule
%		numerical KL&5.3567$\times 10^{-4}$ & $1.5770\times 10^{-5}$ & $4.3684\times 10^{-7}$ \\
		\bottomrule
	\end{tabular}%
	\centering
	\caption{\scriptsize relative error}
	\label{Lplc-table2}%
\end{table}%

%3

\begin{table}[H]
\scriptsize
	\begin{tabular}{l|l|l|l}
		\toprule
		SPOAFD&1.4236$\times 10^{-4}$& 6.5114$\times 10^{-7}$ & $6.0066\times 10^{-8}$\\
		%\midrule
%		numerical KL& 4.6742$\times 10^{-4}$& $9.5706\times 10^{-6}$ & $1.5037\times 10^{-7}$ \\
		\bottomrule
	\end{tabular}%
	\centering
	\caption{\scriptsize relative error}
	\label{Lplc-table3}%
\end{table}%

\begin{figure}[H]
	\begin{minipage}[c]{0.3\textwidth}
		\centering
		\includegraphics[height=4.5cm,width=5cm]{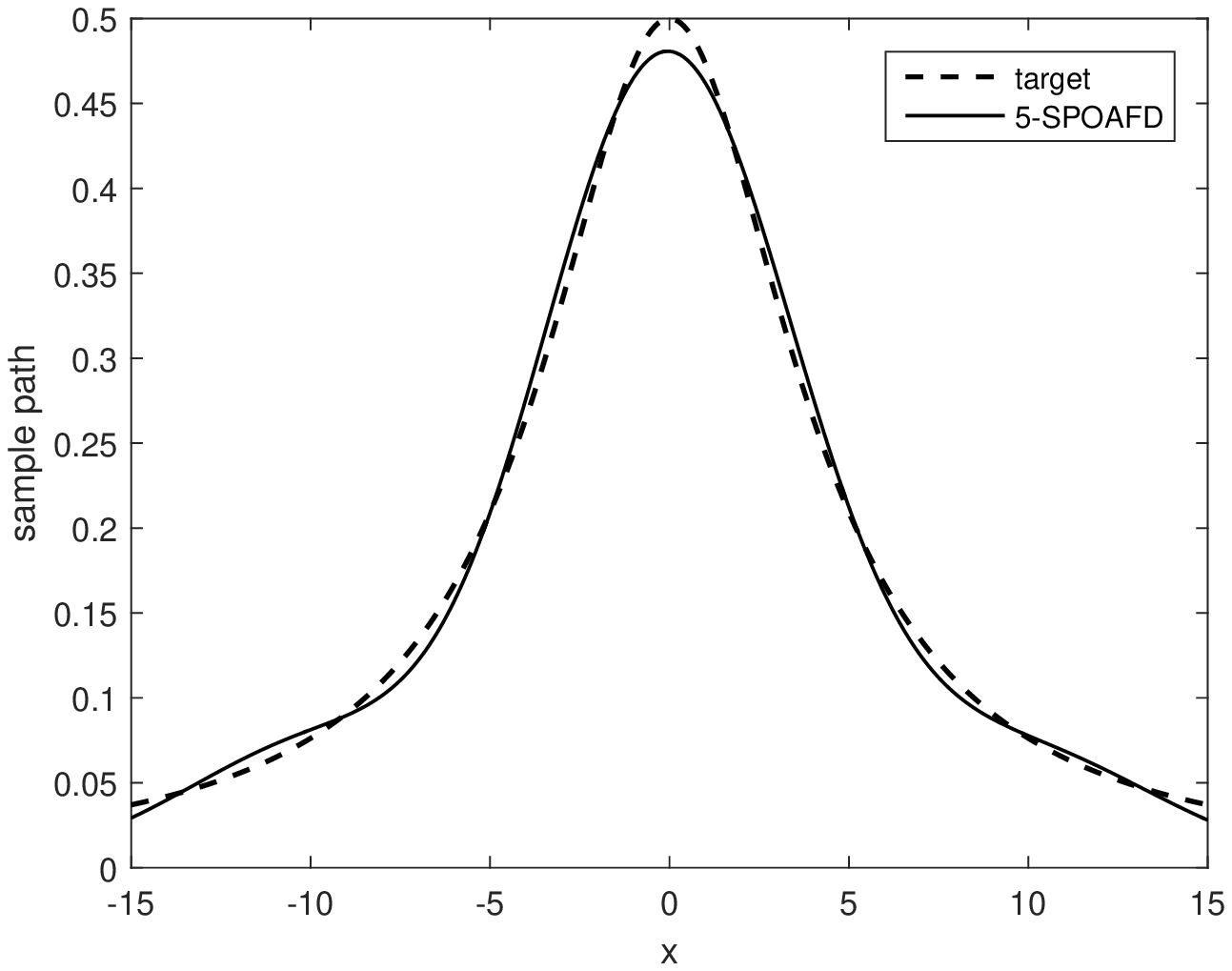}
		\centerline{{\scriptsize 5 partial sum}}
	\end{minipage}
	\begin{minipage}[c]{0.3\textwidth}
		\centering
		\includegraphics[height=4.5cm,width=5cm]{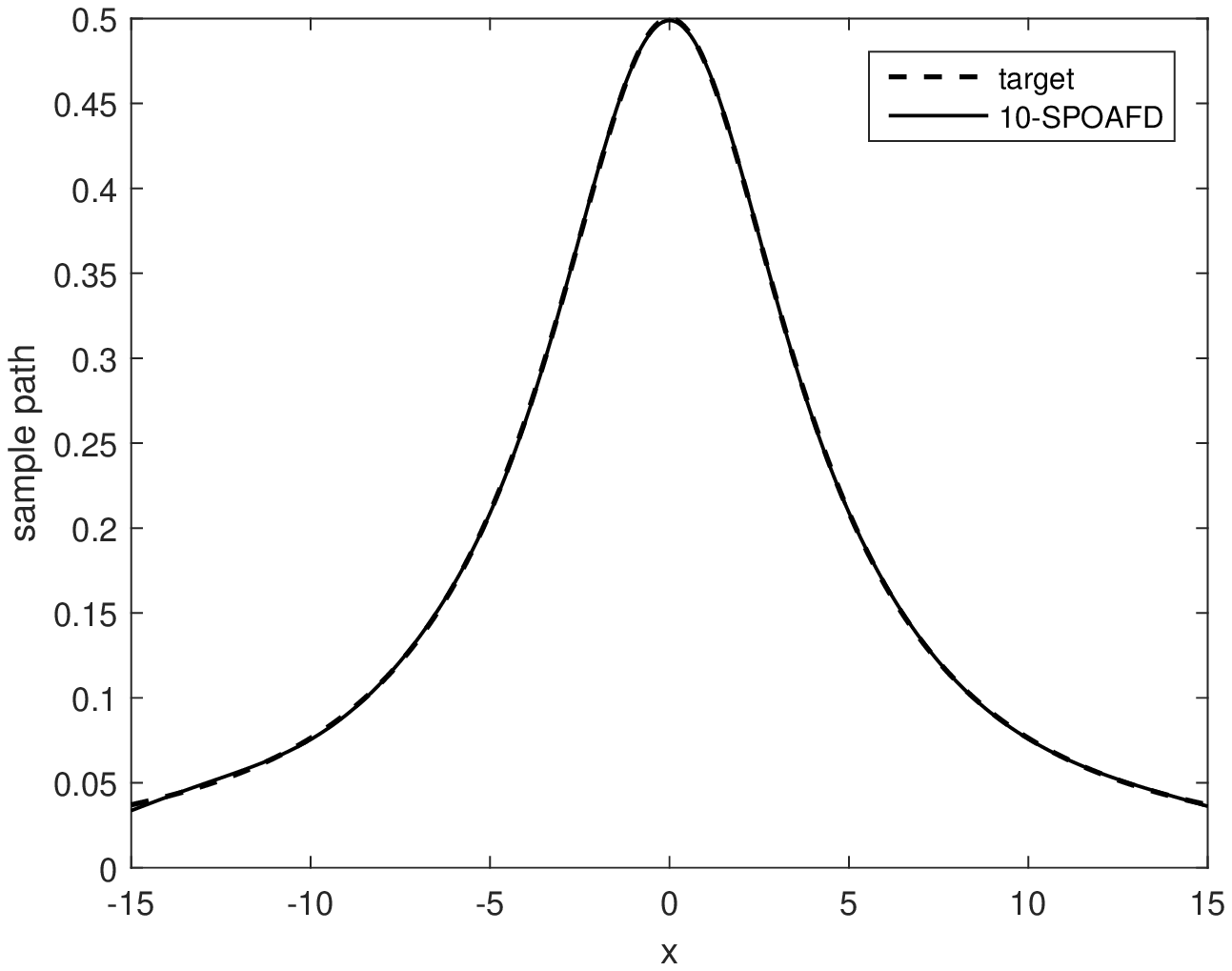}
		\centerline{{\scriptsize 10 partial sum}}
	\end{minipage}
	\begin{minipage}[c]{0.3\textwidth}
		\centering
		\includegraphics[height=4.5cm,width=5cm]{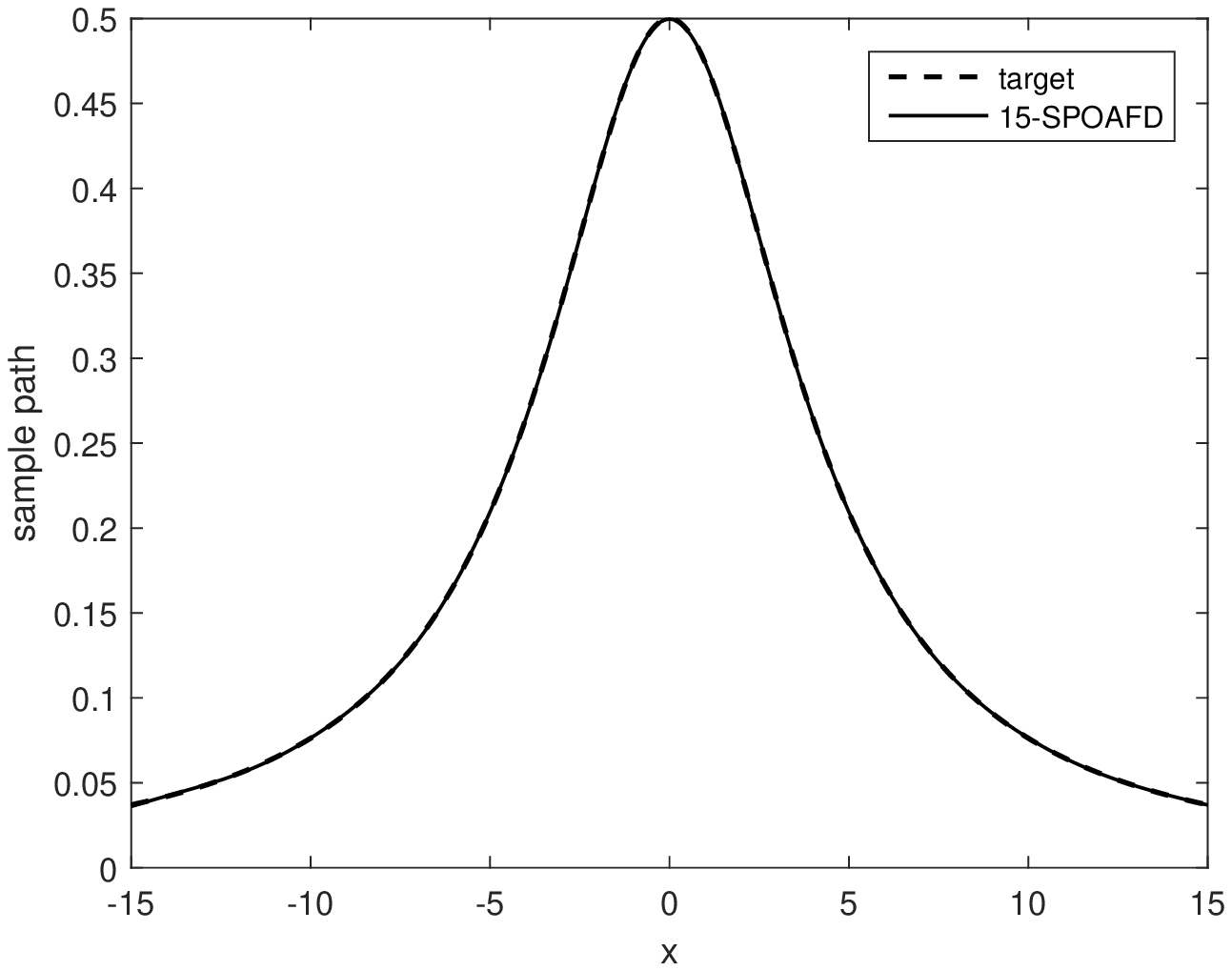}
		\centerline{{\scriptsize 15 partial sum}}
	\end{minipage}
	
	\caption{$X=0$}
\label{figure4}
\end{figure}

\begin{figure}[H]
	\begin{minipage}[c]{0.3\textwidth}
		\centering
		\includegraphics[height=4.5cm,width=5cm]{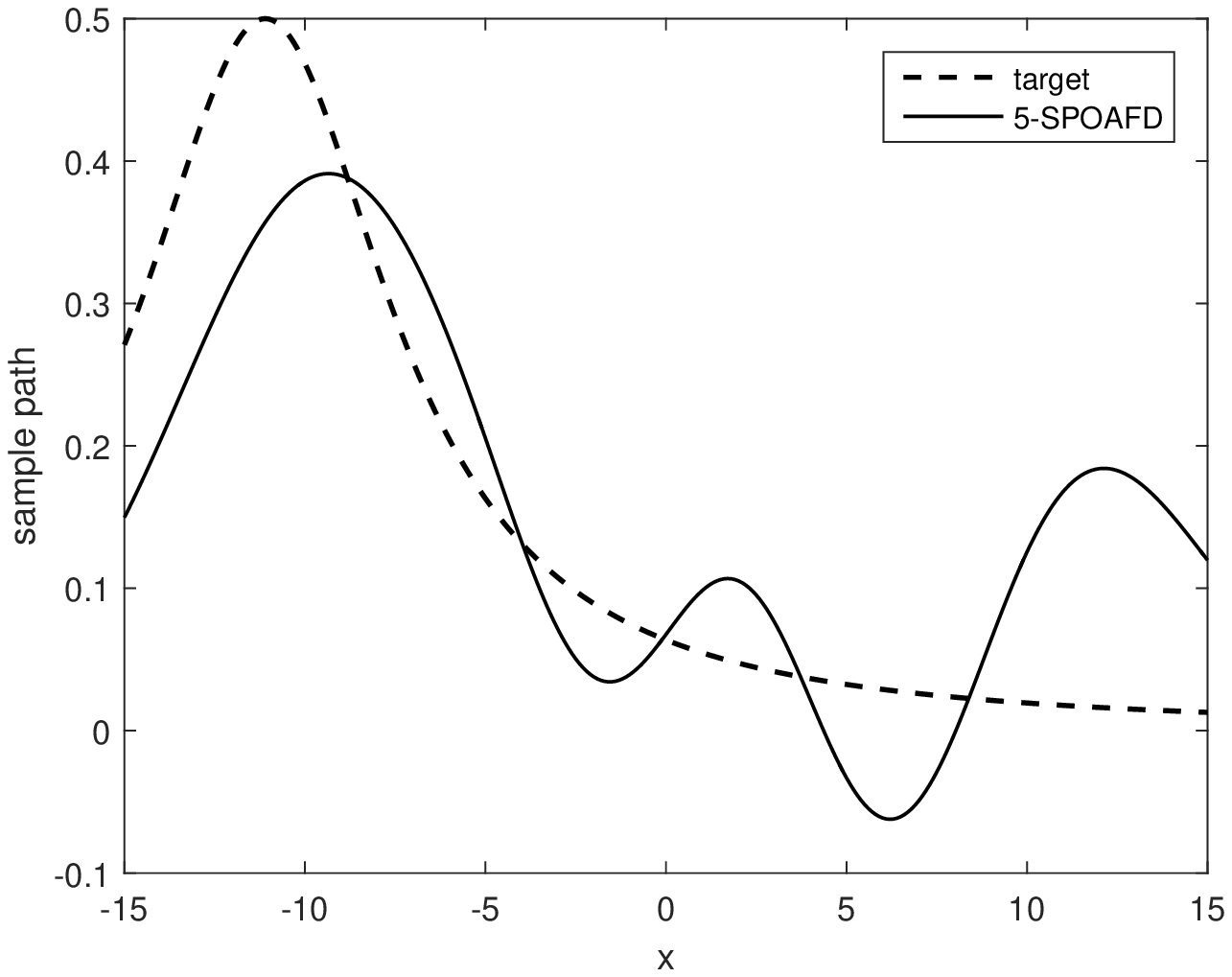}
		\centerline{{\scriptsize 5 partial sum}}
	\end{minipage}
	\begin{minipage}[c]{0.3\textwidth}
		\centering
		\includegraphics[height=4.5cm,width=5cm]{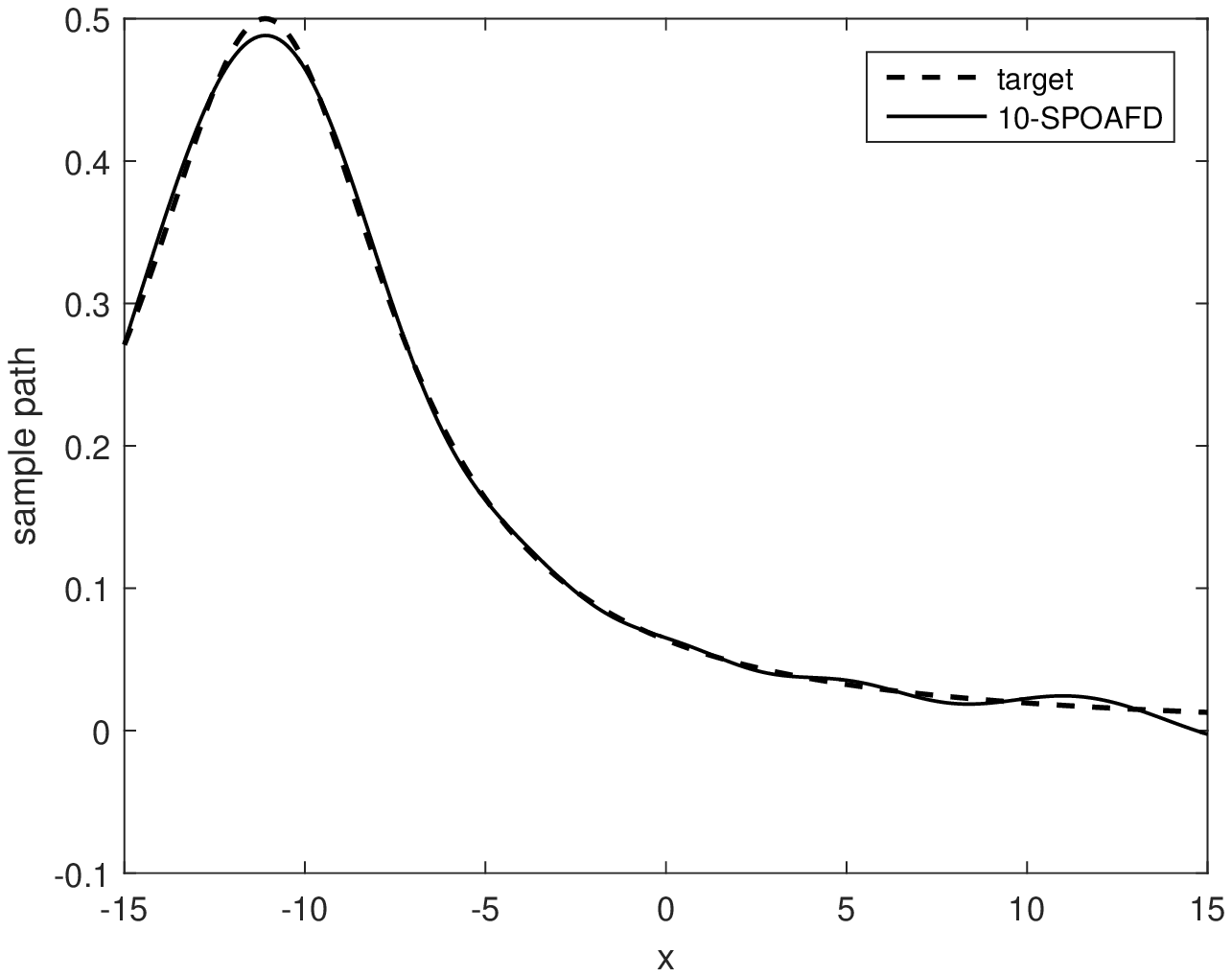}
		\centerline{{\scriptsize  10 partial sum}}
	\end{minipage}
	\begin{minipage}[c]{0.3\textwidth}
		\centering
		\includegraphics[height=4.5cm,width=5cm]{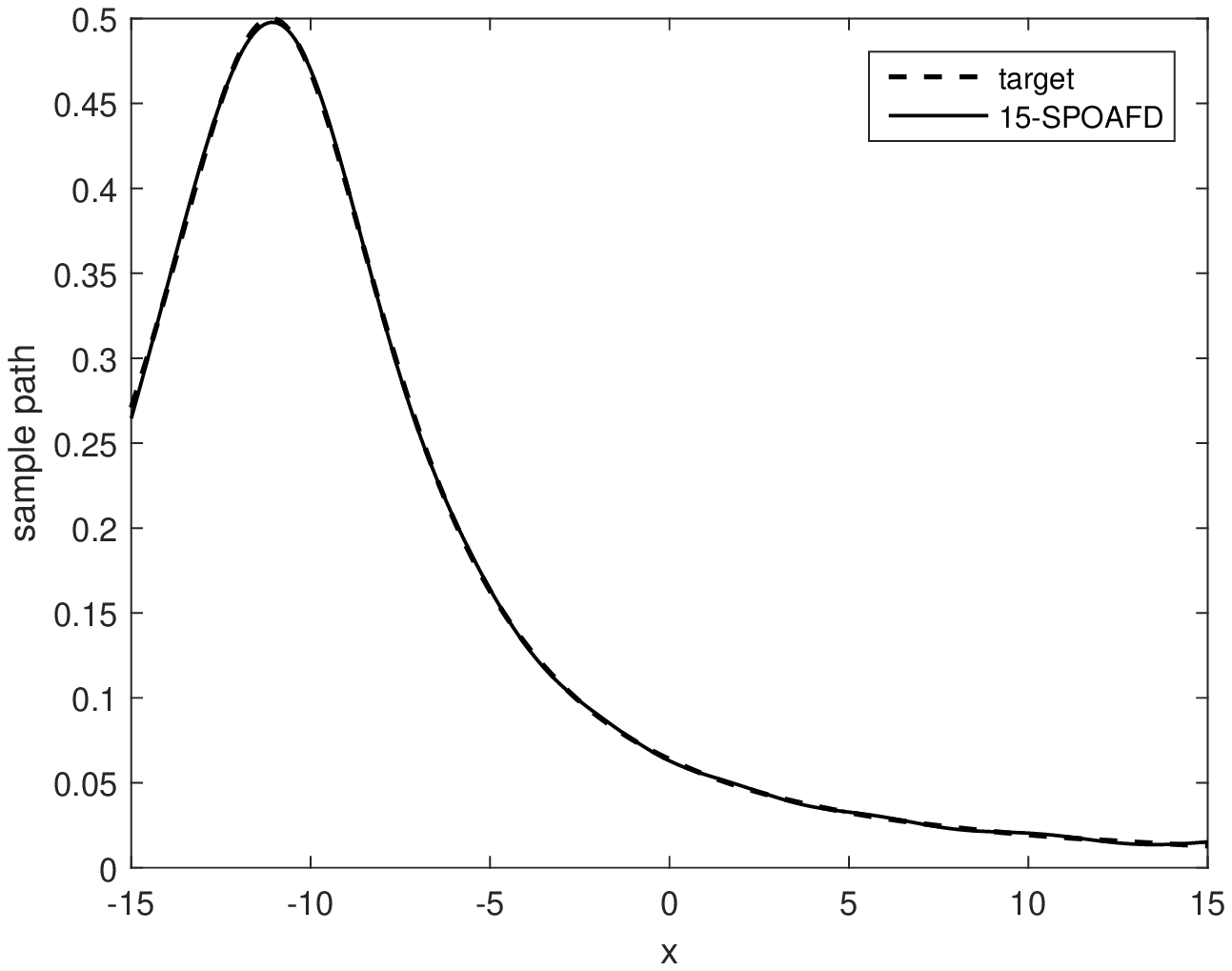}
		\centerline{{\scriptsize  15 partial sum}}
	\end{minipage}
	
	\caption{$X=-3.7$}
\label{figure5}
\end{figure}

\begin{figure}[H]
	\begin{minipage}[c]{0.3\textwidth}
		\centering
		\includegraphics[height=4.5cm,width=5cm]{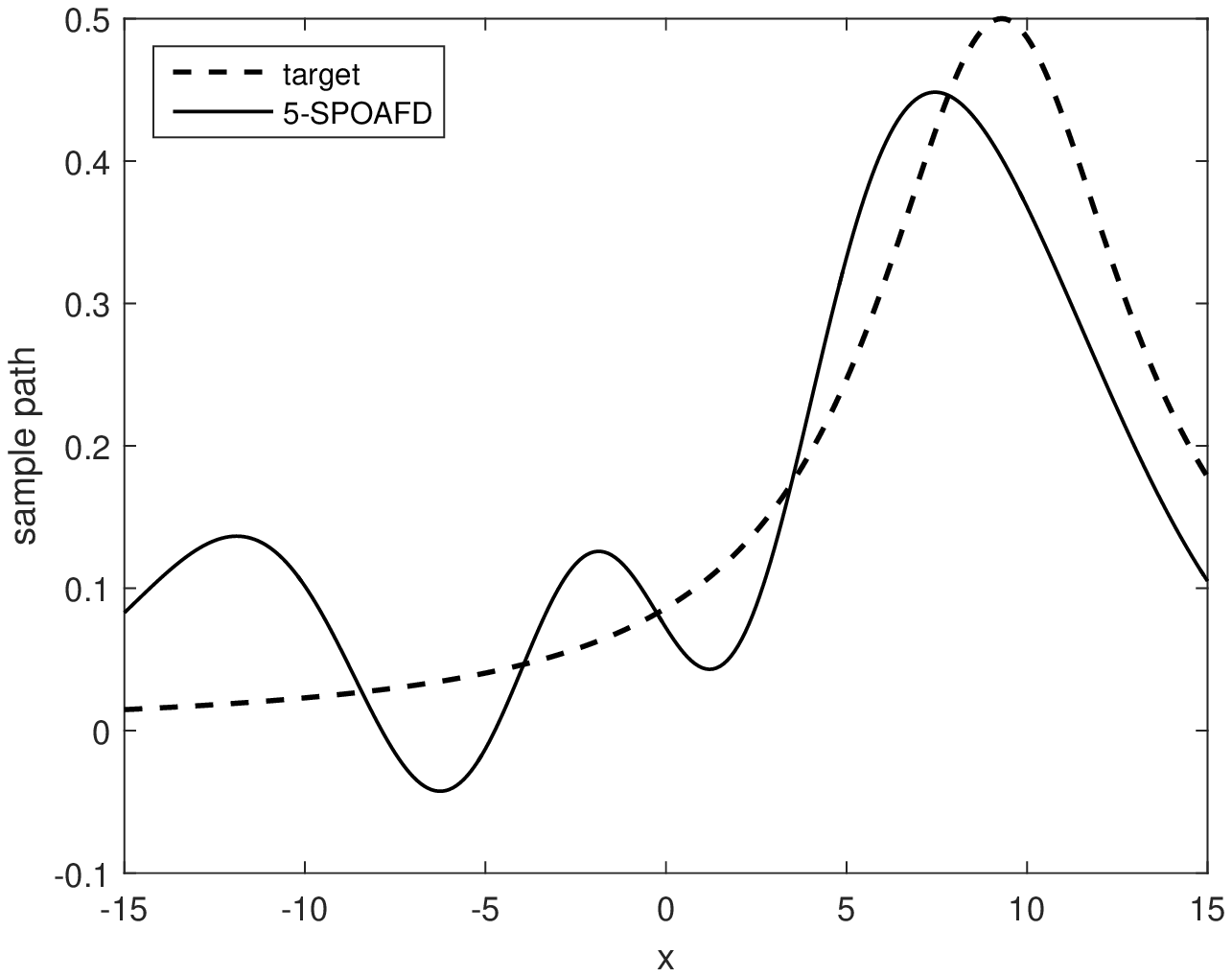}
		\centerline{{\scriptsize  5 partial sum}}
	\end{minipage}
	\begin{minipage}[c]{0.3\textwidth}
		\centering
		\includegraphics[height=4.5cm,width=5cm]{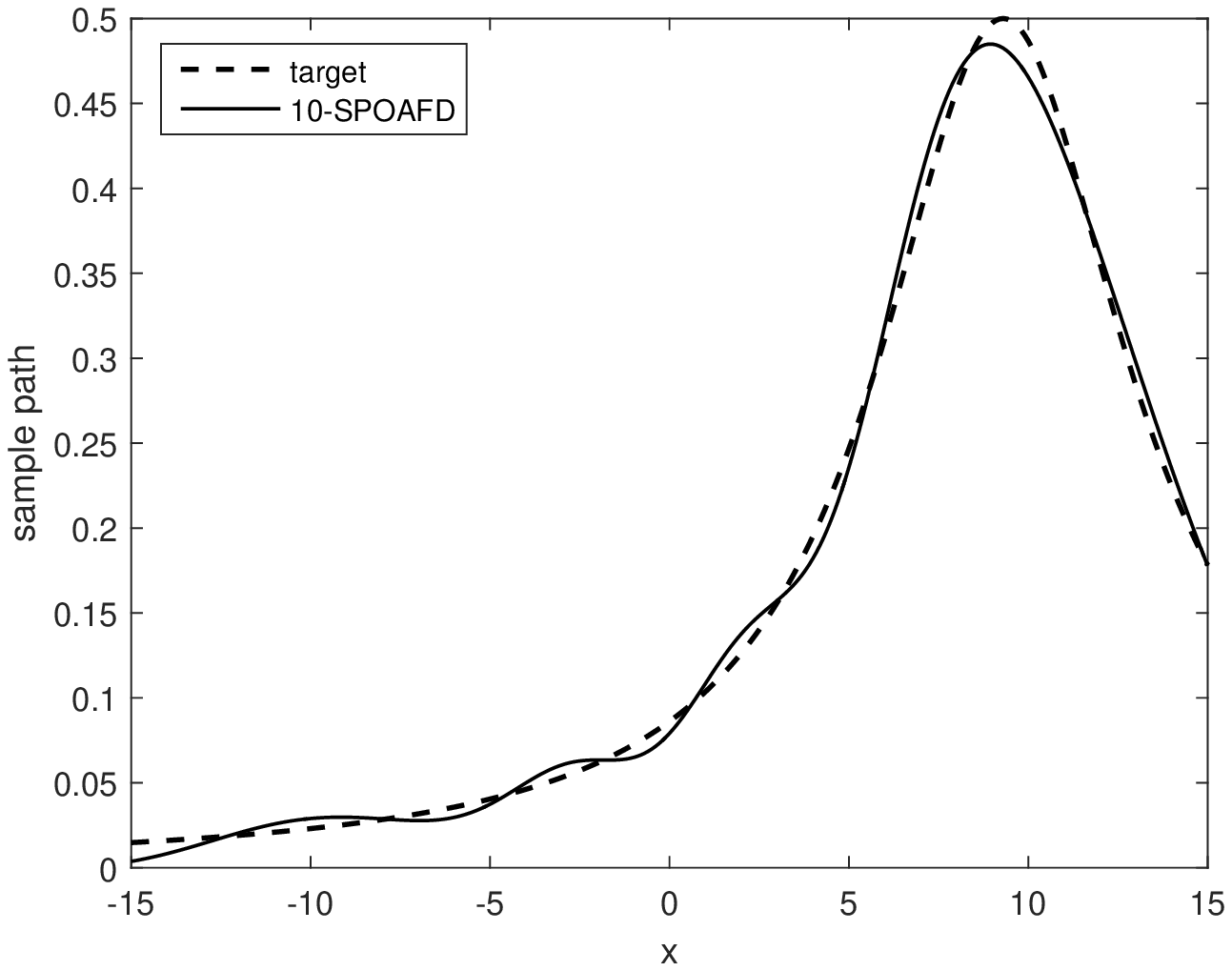}
		\centerline{{\scriptsize  10 partial sum}}
	\end{minipage}
	\begin{minipage}[c]{0.3\textwidth}
		\centering
		\includegraphics[height=4.5cm,width=5cm]{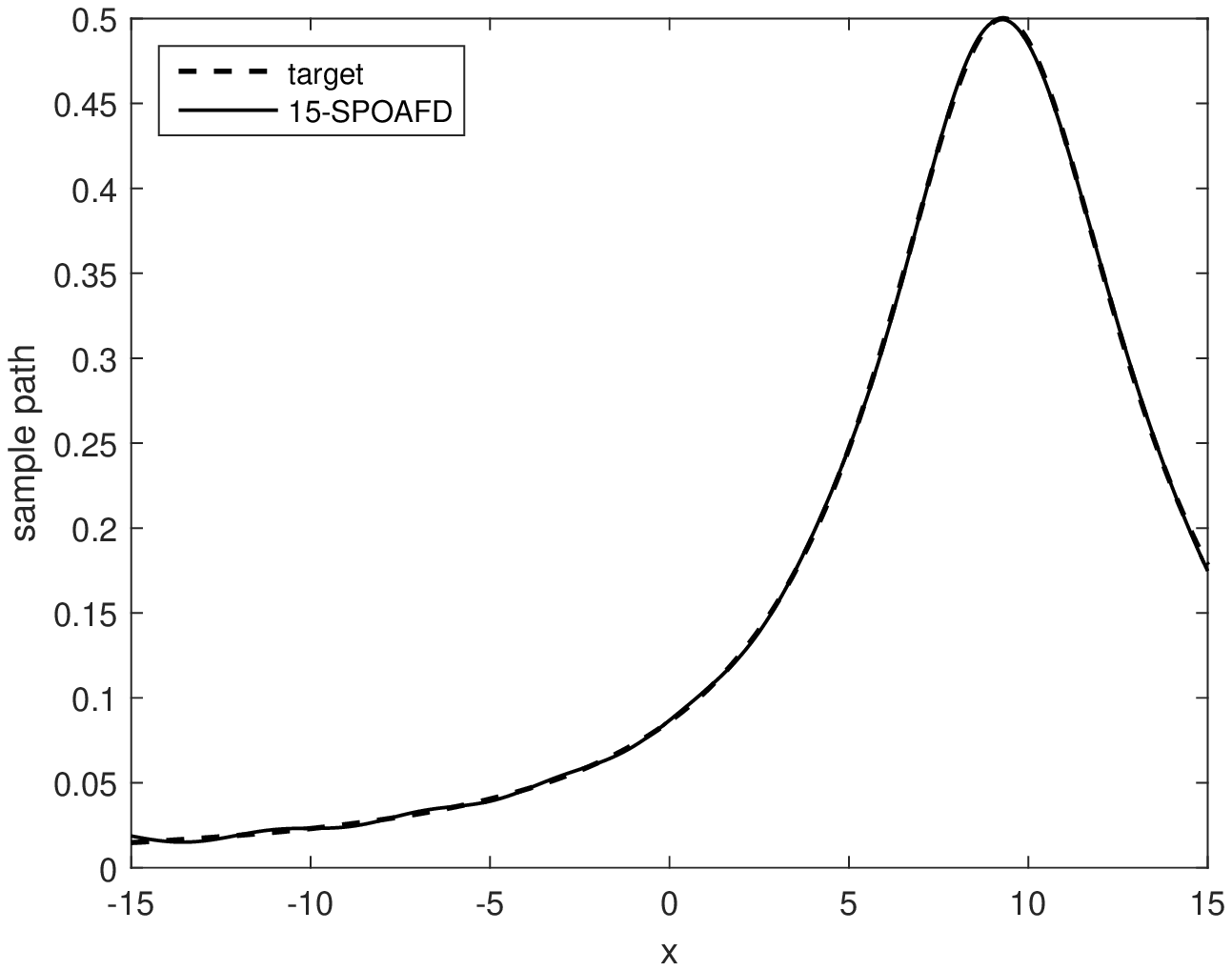}
		\centerline{{\scriptsize  15 partial sum}}
	\end{minipage}
	\caption{$X=3.1$}
\label{figure6}
\end{figure}

\begin{Exmp}\label{Exmp-heat}
The other case is the heat equation with the random initial condition:
\begin{equation}
\left\{
\begin{aligned}
&(\partial_t-\Delta_x )u(t,x,\omega) =0, ~ ~(t,x) \in \mathbb{R}^{1}_+\times \mathbb{R}^1,~{\rm a.s.}~\omega \in \Omega,\\
&u(0,x,\omega)=\frac{1}{2+(\frac{x}{3}-X)^2}},~ ~x \in \mathbb{R}^1,~{\rm a.s.},
{\end{aligned}
\right.
\end{equation}
where $X$ is the evaluation of the random variable with the density function
$p(s)=\frac{1}{\sqrt{2\pi}}e^{-\frac{s^2}{2}}.$
\end{Exmp}
 Again for this example we can make direct use of the distribution of $X.$ We have three groups of figures respectively for visual effects of the experiments at those $\omega$ for which $X=0, -3.7,$ and $3.1$, which can be shown in Figure \ref{figure4}, \ref{figure5} and \ref{figure6}. The relative errors in Table \ref{heat-table1},\ref{heat-table2},\ref{heat-table3} are computed according to the formulas given in the relevant algorithms.

%1
%\begin{comment}

\begin{table}[H]
\scriptsize
	\begin{tabular}{l|l|l|l}
		\toprule
		SPOAFD&0.0014 & $9.0613\times 10^{-6}$ & $7.1549\times 10^{-7}$ \\
		%\midrule
%		Hermite expansion&0.0015 & $1.5138\times 10^{-4}$ & $5.1368\times 10^{-6}$ \\
		\bottomrule
	\end{tabular}%
	\centering
	\caption{\scriptsize relative error}
	\label{heat-table1}%
\end{table}%

%2
%\end{comment}

\begin{table}[H]
\scriptsize
	\begin{tabular}{l|l|l|l}
		\toprule
		SPOAFD&0.1759 & 4.8514$\times 10^{-4}$ & $2.7590\times 10^{-5}$\\
		%\midrule
%		Hermite expansion&1.3774 &  0.0528 & 0.0014 \\
		\bottomrule
	\end{tabular}%
	\centering
	\caption{\scriptsize relative error}
	\label{heat-table2}%
\end{table}%

%3
%\begin{comment}
%\end{comment}

\begin{table}[H]
\scriptsize
	\begin{tabular}{l|l|l|l}
		\toprule
		SPOAFD&0.1064 & 0.0016& $2.5211\times 10^{-5}$\\
		%\midrule
%		Hermite expansion&0.0950 &  0.0048 & $4.2147\times 10^{-4}$ \\
		\bottomrule
	\end{tabular}%
	\centering
	\caption{\scriptsize relative error}
	\label{heat-table3}%
\end{table}%

\begin{figure}[H]
	\begin{minipage}[c]{0.3\textwidth}
		\centering
		\includegraphics[height=4.5cm,width=5cm]{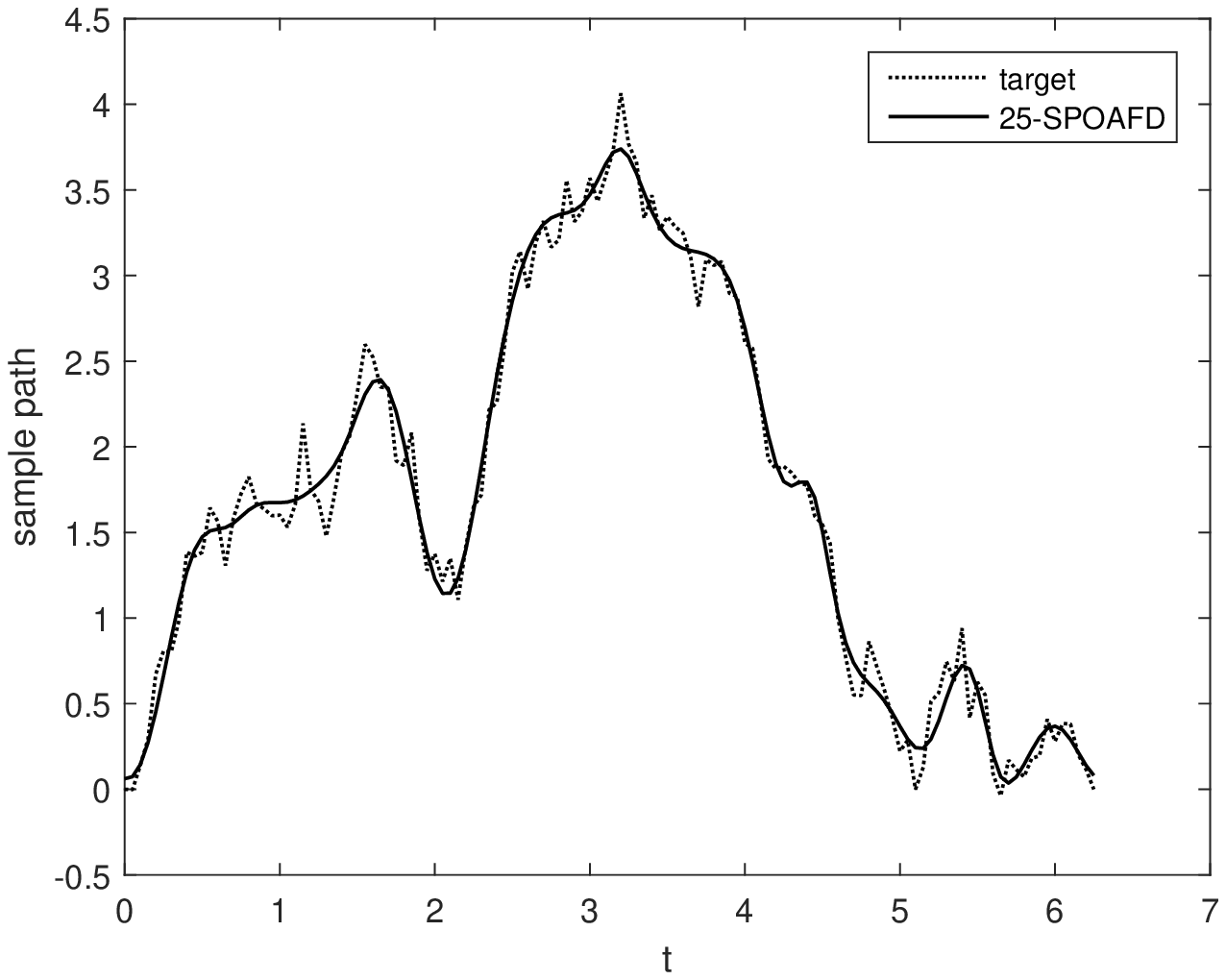}
		\centerline{{\scriptsize SPOAFD: 25 partial sum}}
	\end{minipage}
	\begin{minipage}[c]{0.3\textwidth}
		\centering
		\includegraphics[height=4.5cm,width=5cm]{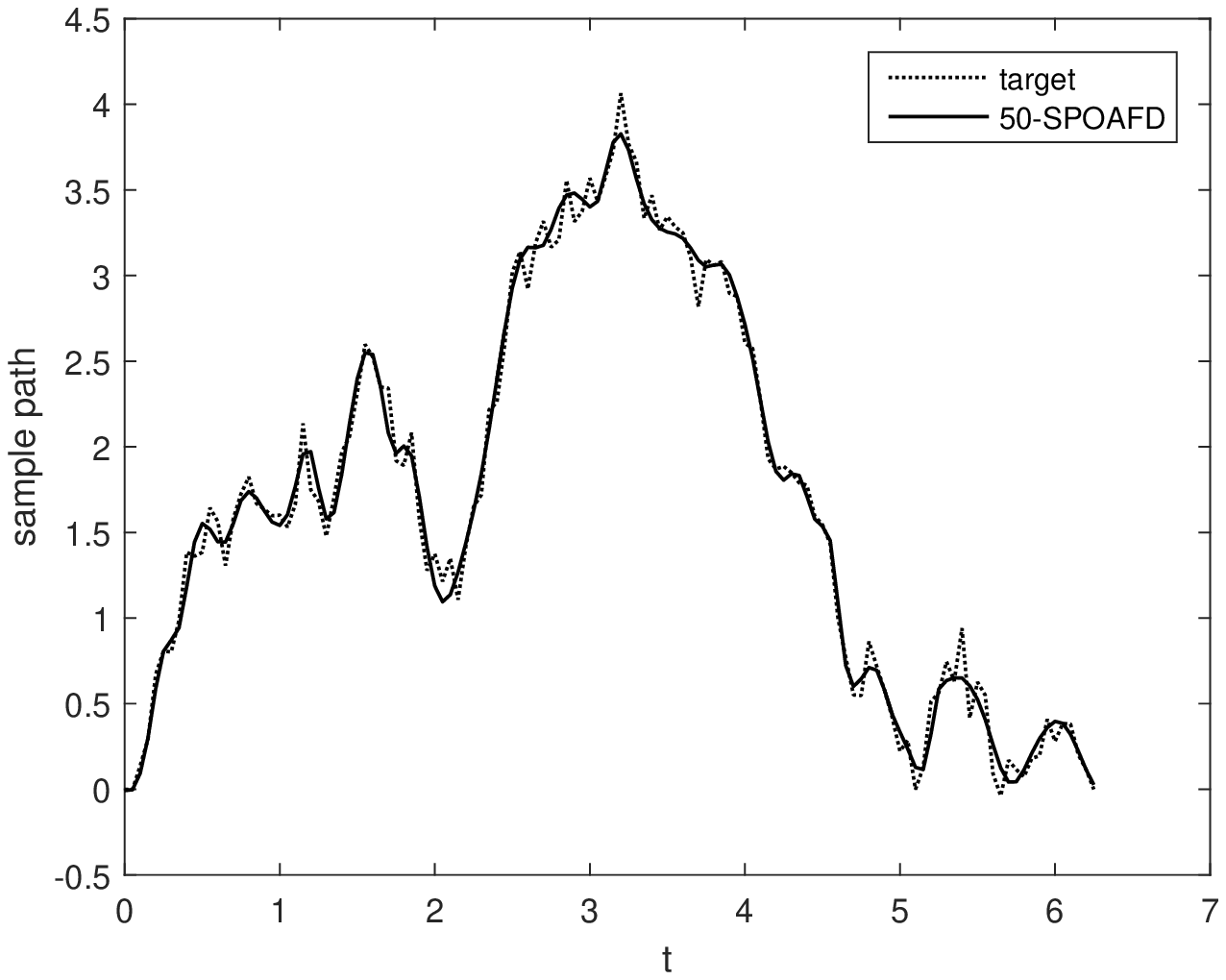}
		\centerline{{\scriptsize SPOAFD: 50 partial sum}}
	\end{minipage}
	\begin{minipage}[c]{0.3\textwidth}
		\centering
		\includegraphics[height=4.5cm,width=5cm]{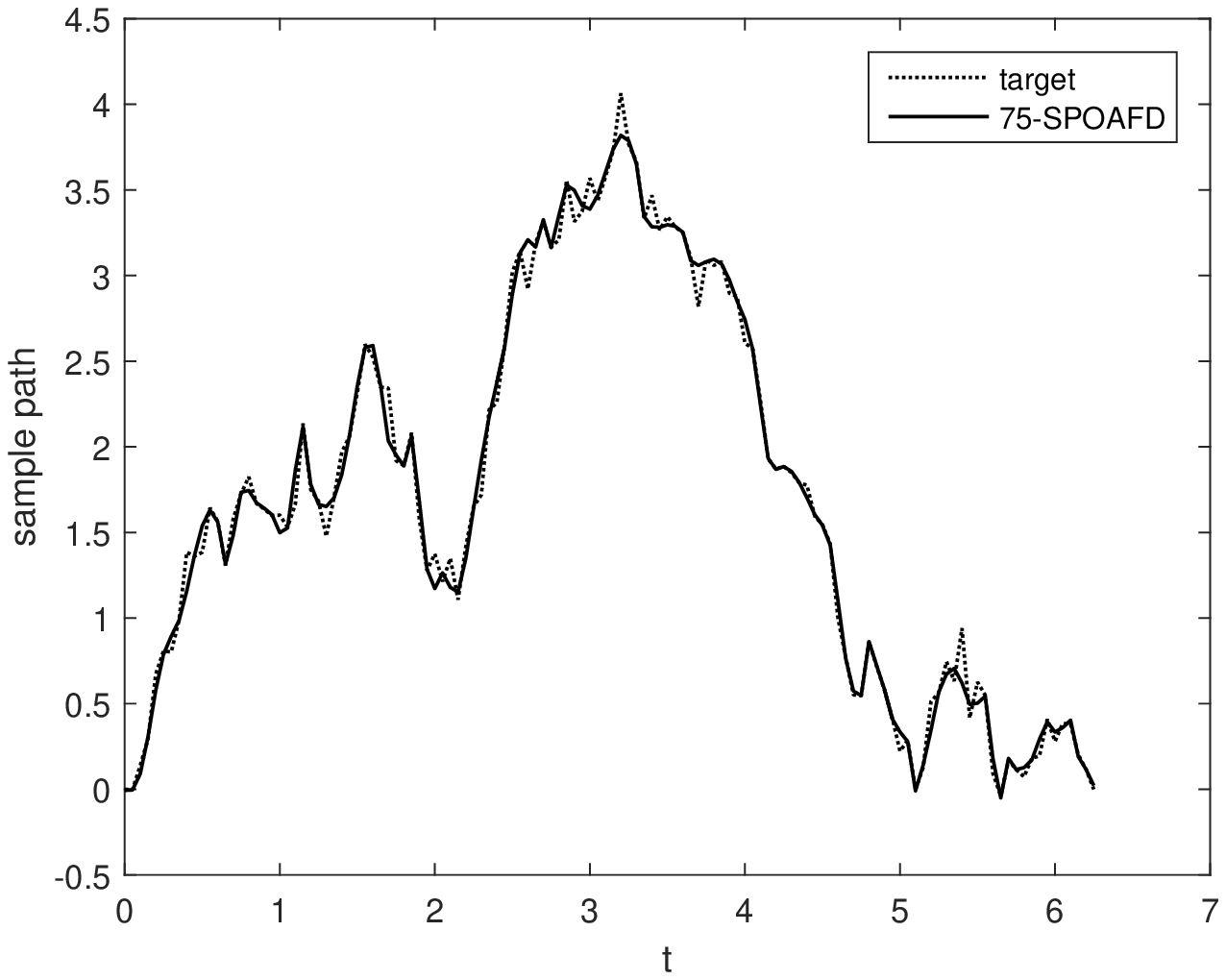}
		\centerline{{\scriptsize SPOAFD: 75 partial sum}}
	\end{minipage}
	
	\begin{minipage}[c]{0.3\textwidth}
		\centering
		\includegraphics[height=4.5cm,width=5cm]{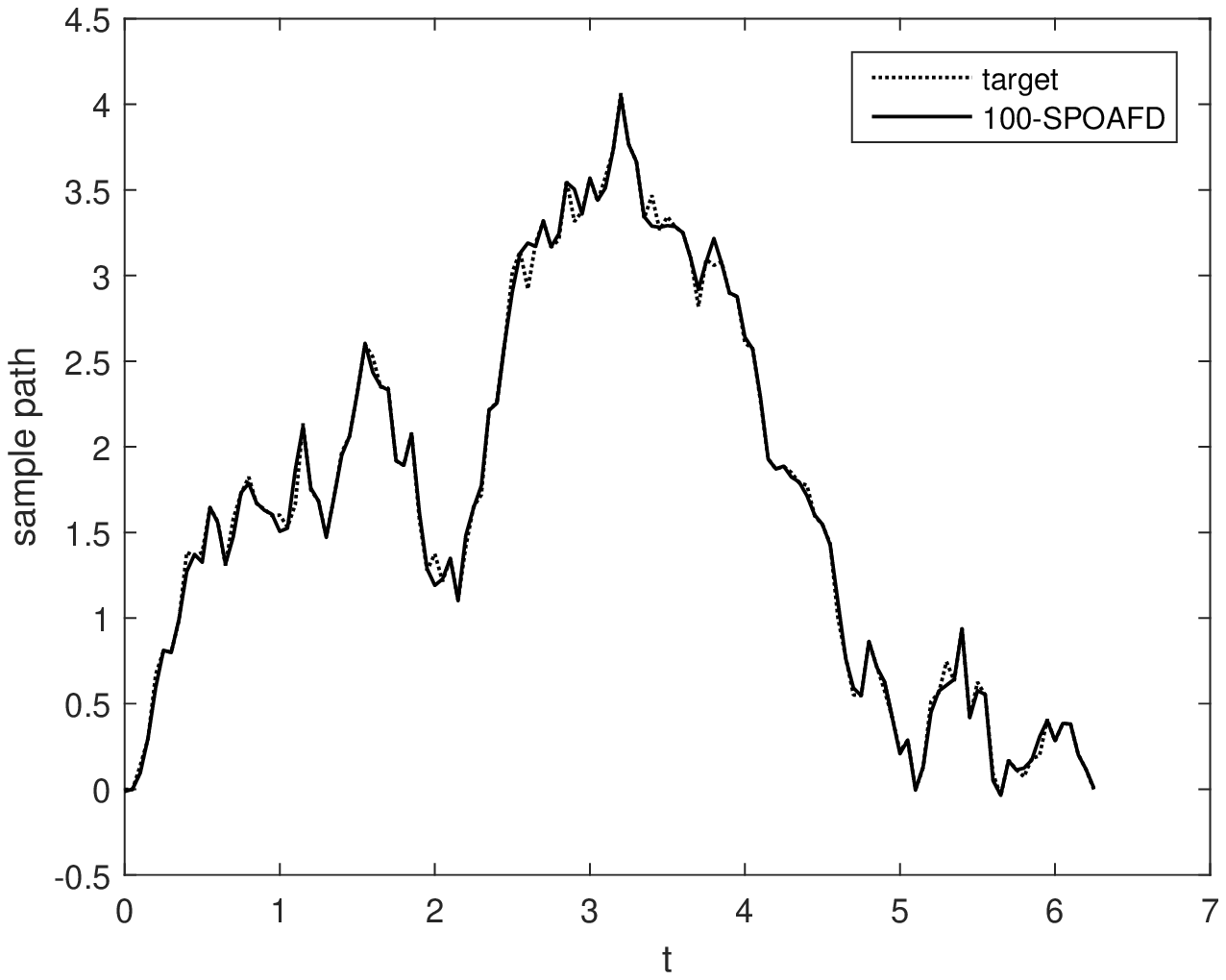}
		\centerline{{\scriptsize SPOAFD: 100 partial sum}}
	\end{minipage}
	\begin{minipage}[c]{0.3\textwidth}
		\centering
		\includegraphics[height=4.5cm,width=5cm]{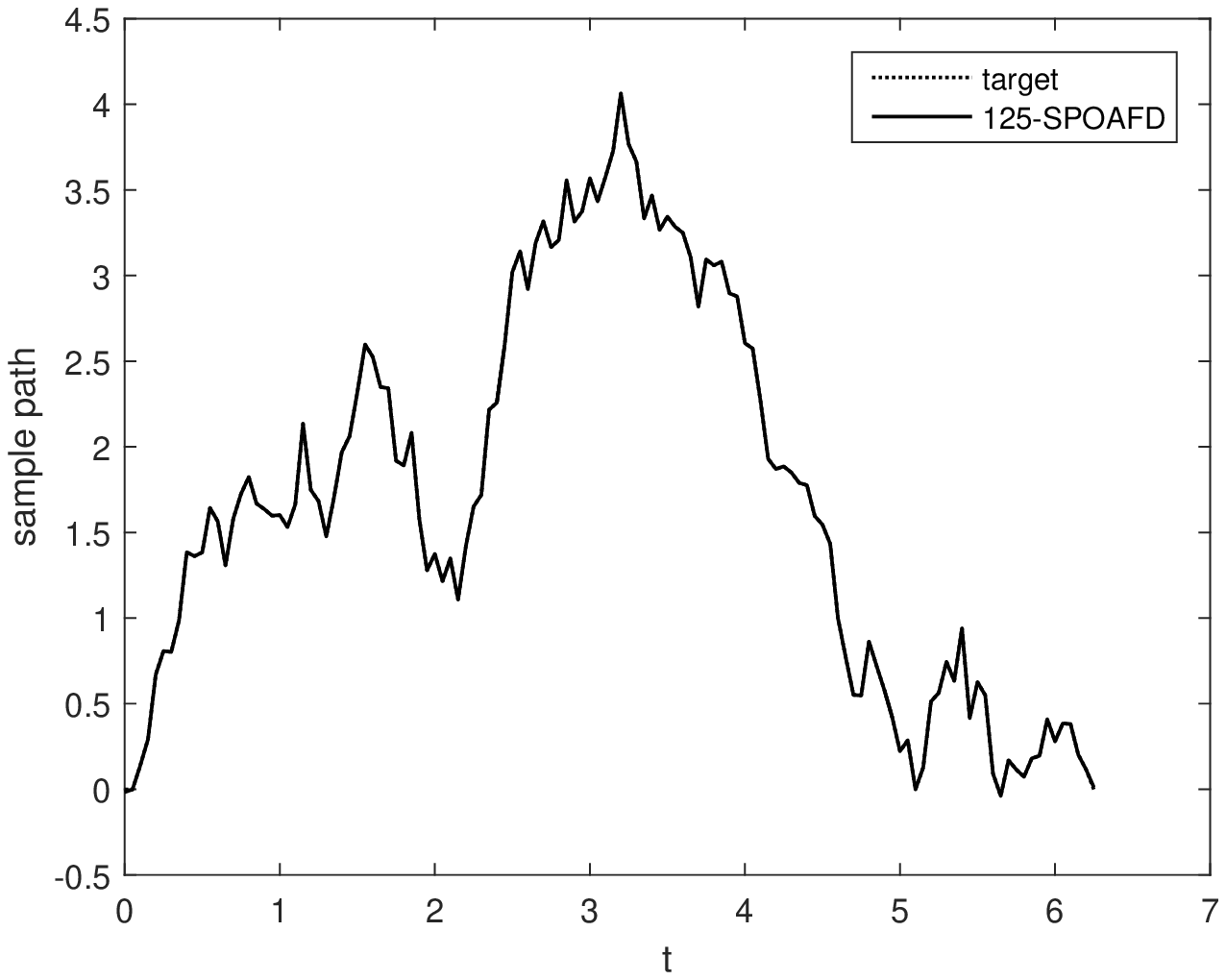}
		\centerline{{\scriptsize SPOAFD2: 125 partial sum}}
	\end{minipage}
	\begin{minipage}[c]{0.3\textwidth}
		\centering
		\includegraphics[height=4.5cm,width=5cm]{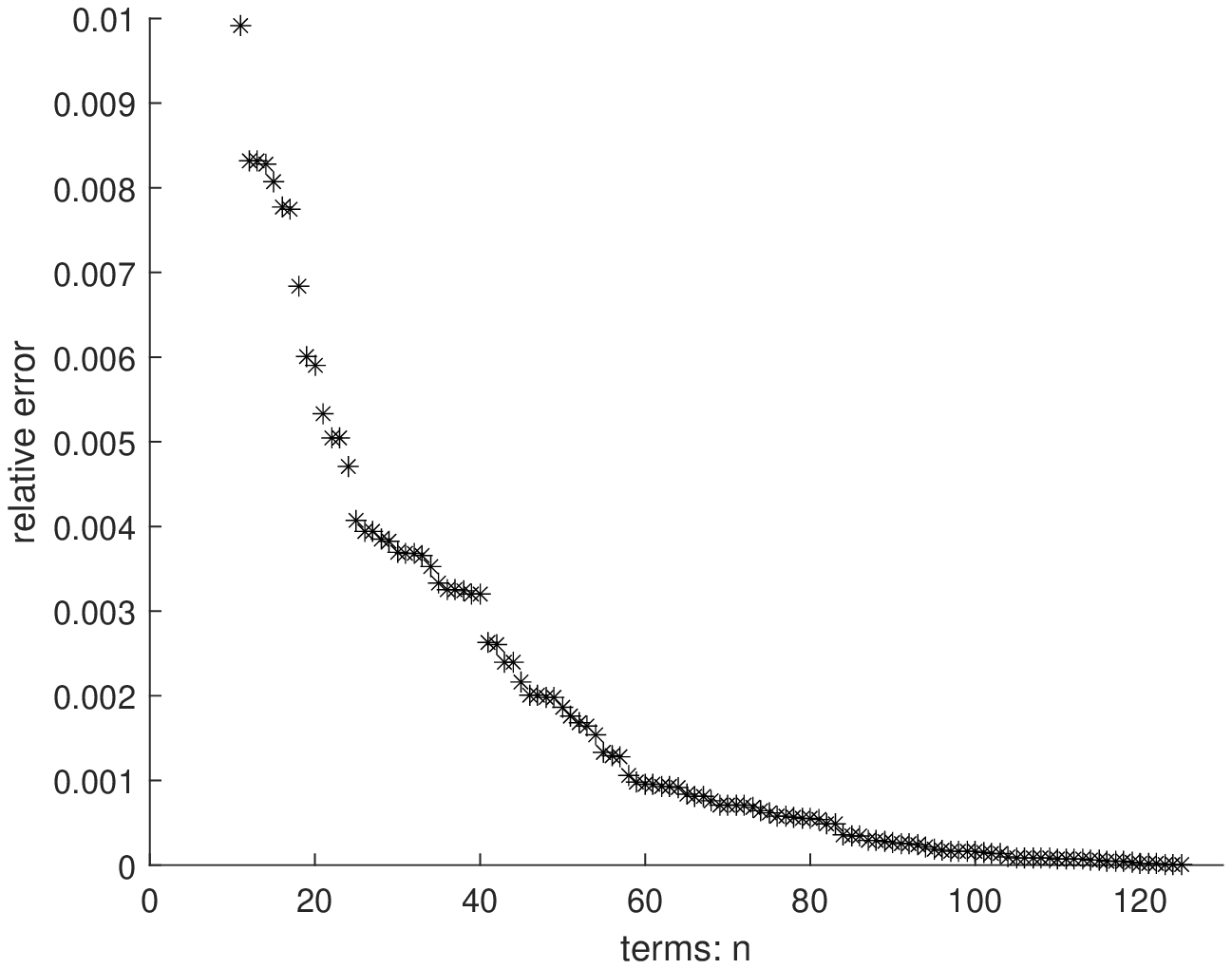}
		\centerline{{\scriptsize relative error}}
	\end{minipage}
	\caption{sample path I of Brownian bridge}
\label{figureB1}
\end{figure}

\begin{figure}[H]
	\begin{minipage}[c]{0.3\textwidth}
		\centering
		\includegraphics[height=4.5cm,width=5cm]{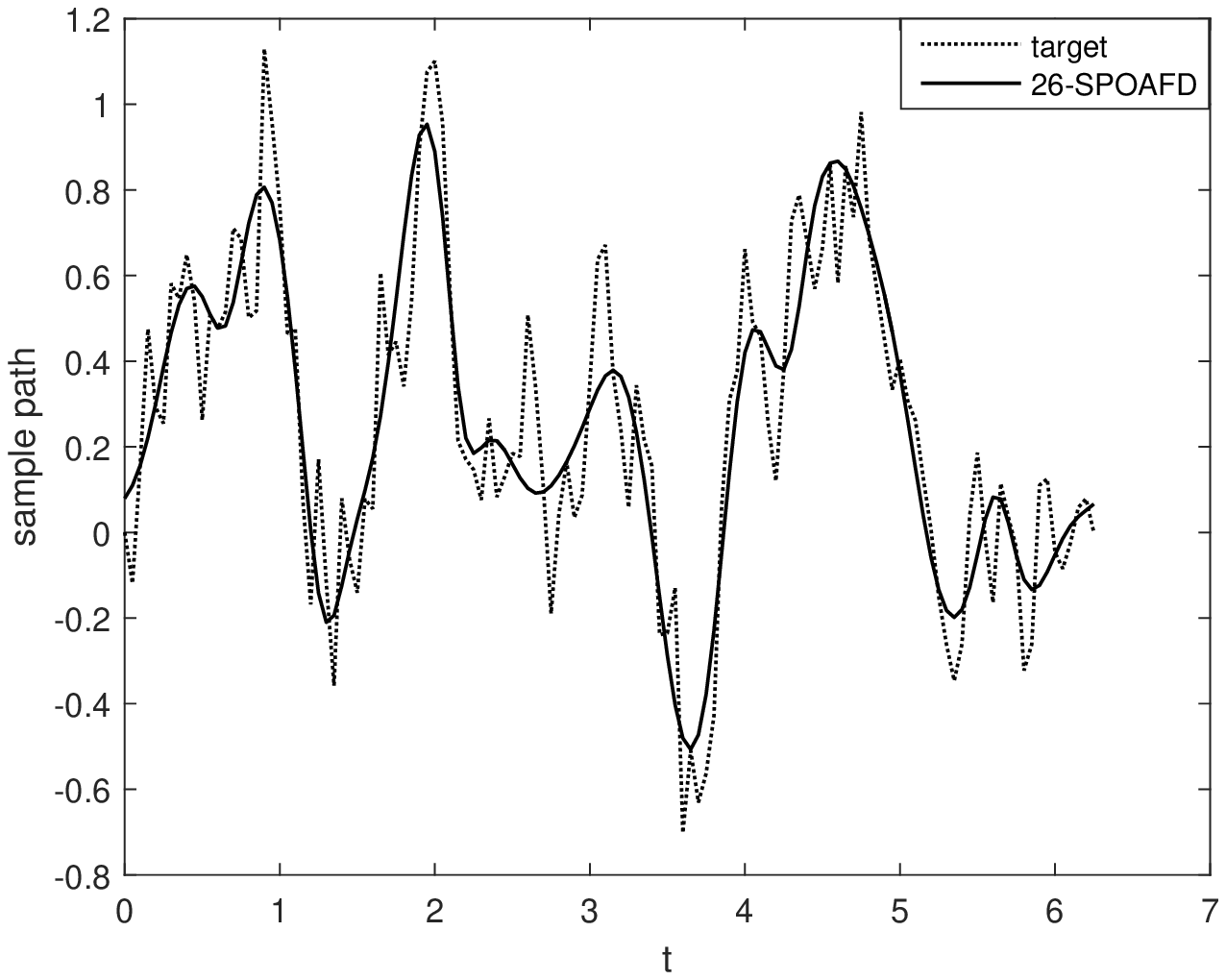}
		\centerline{{\scriptsize SPOAFD: 26 partial sum}}
	\end{minipage}
	\begin{minipage}[c]{0.3\textwidth}
		\centering
		\includegraphics[height=4.5cm,width=5cm]{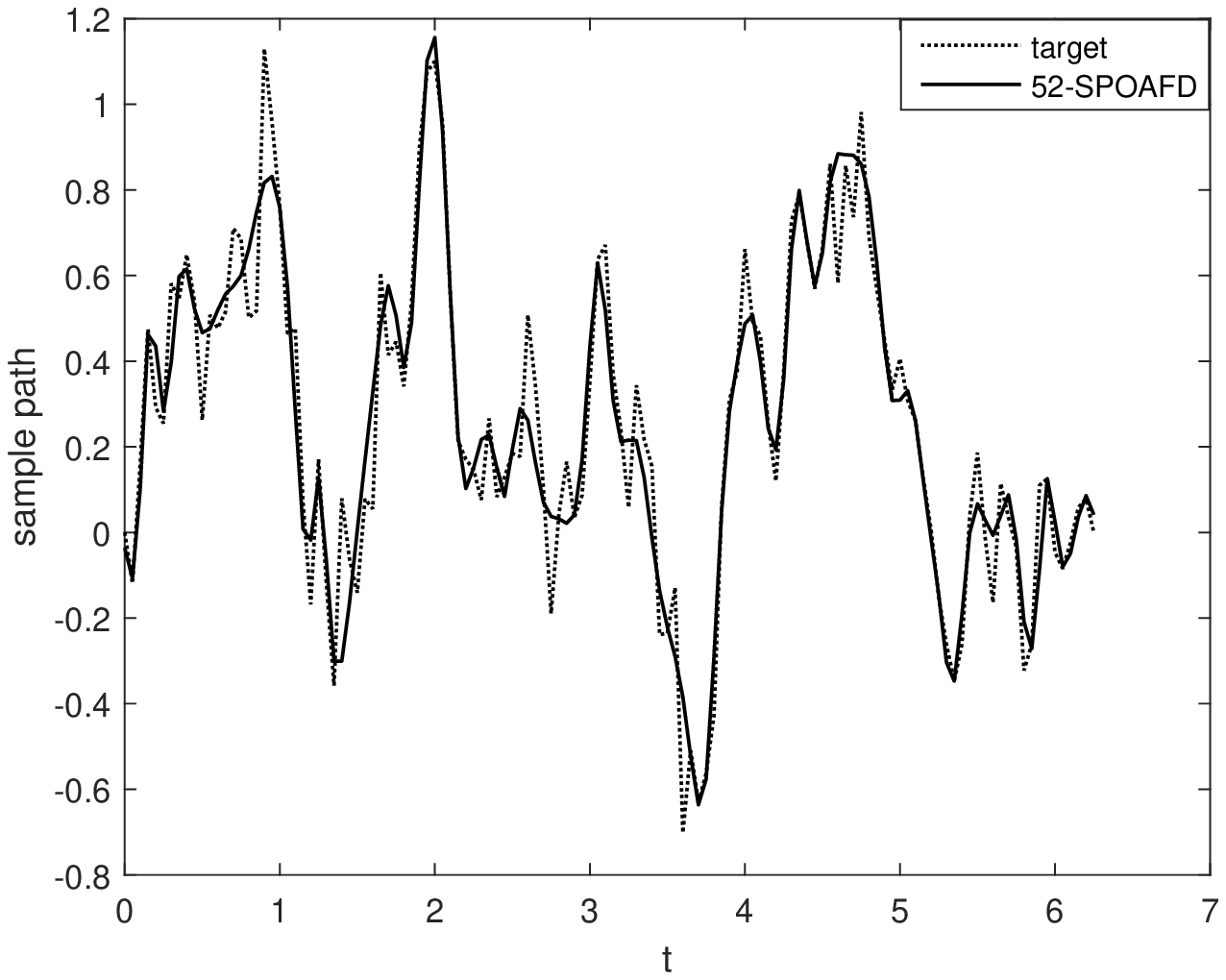}
		\centerline{{\scriptsize SPOAFD: 52 partial sum}}
	\end{minipage}
	\begin{minipage}[c]{0.3\textwidth}
		\centering
		\includegraphics[height=4.5cm,width=5cm]{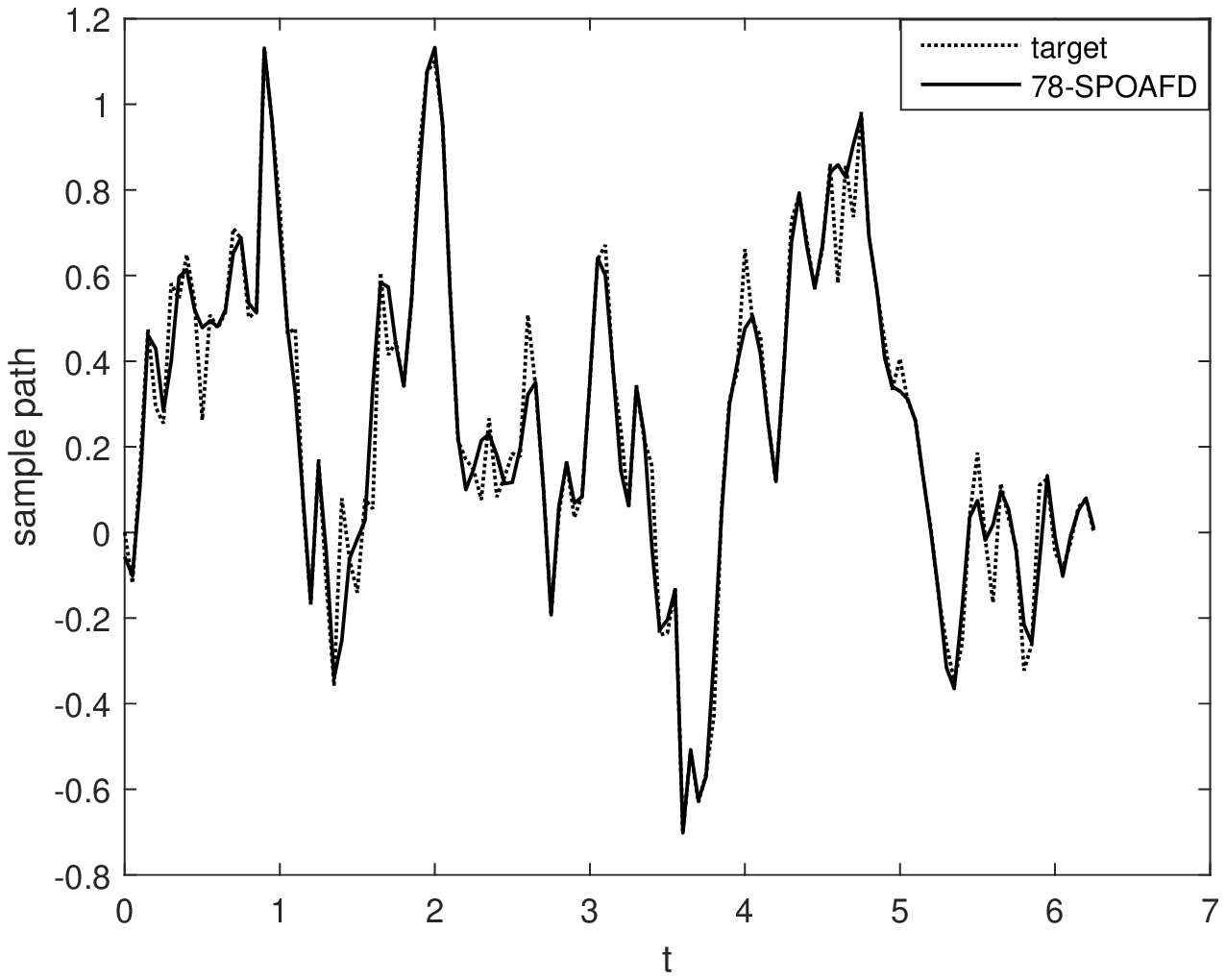}
		\centerline{{\scriptsize SPOAFD: 78 partial sum}}
	\end{minipage}
	
	\begin{minipage}[c]{0.3\textwidth}
		\centering
		\includegraphics[height=4.5cm,width=5cm]{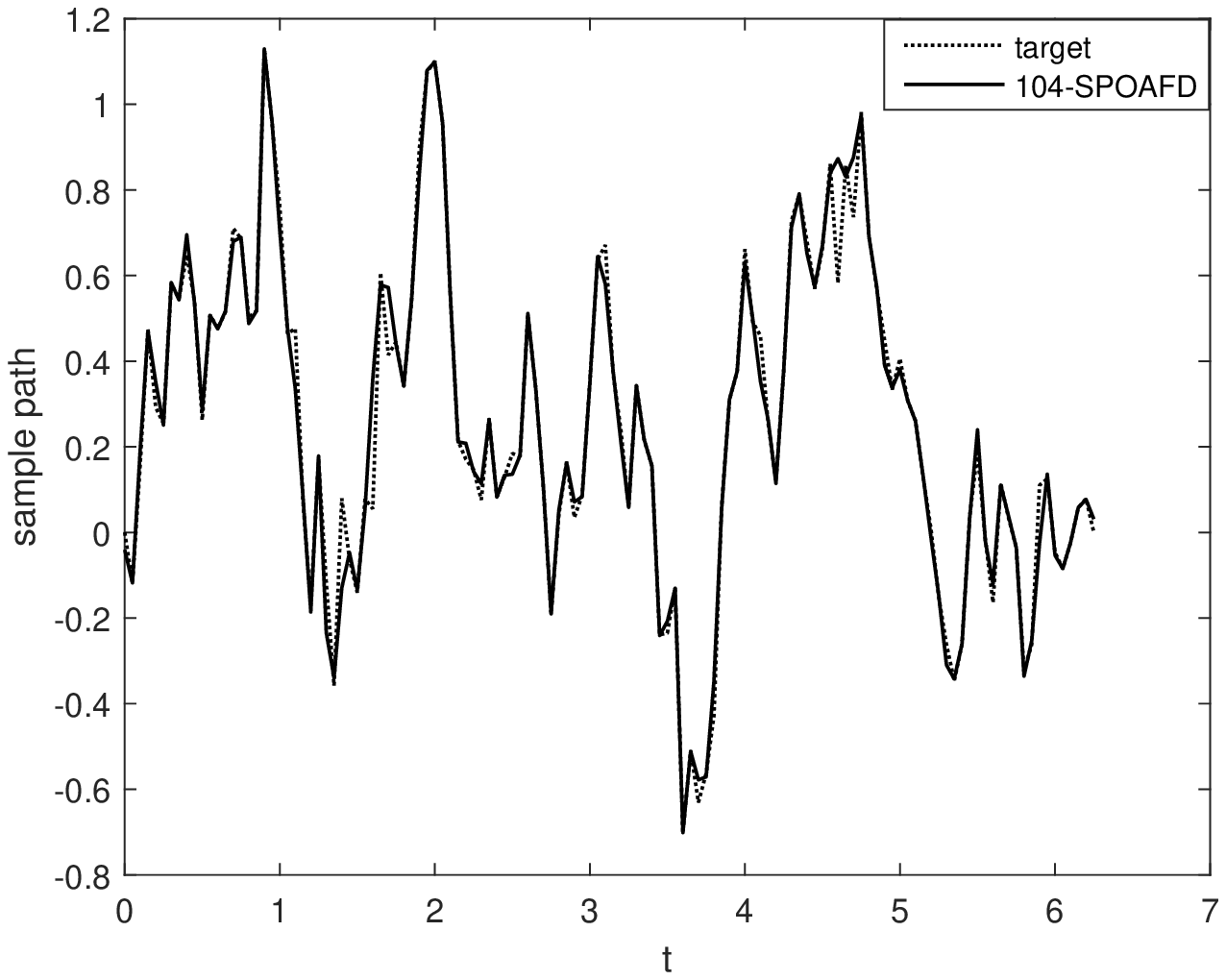}
		\centerline{{\scriptsize SPOAFD: 104 partial sum}}
	\end{minipage}
	\begin{minipage}[c]{0.3\textwidth}
		\centering
		\includegraphics[height=4.5cm,width=5cm]{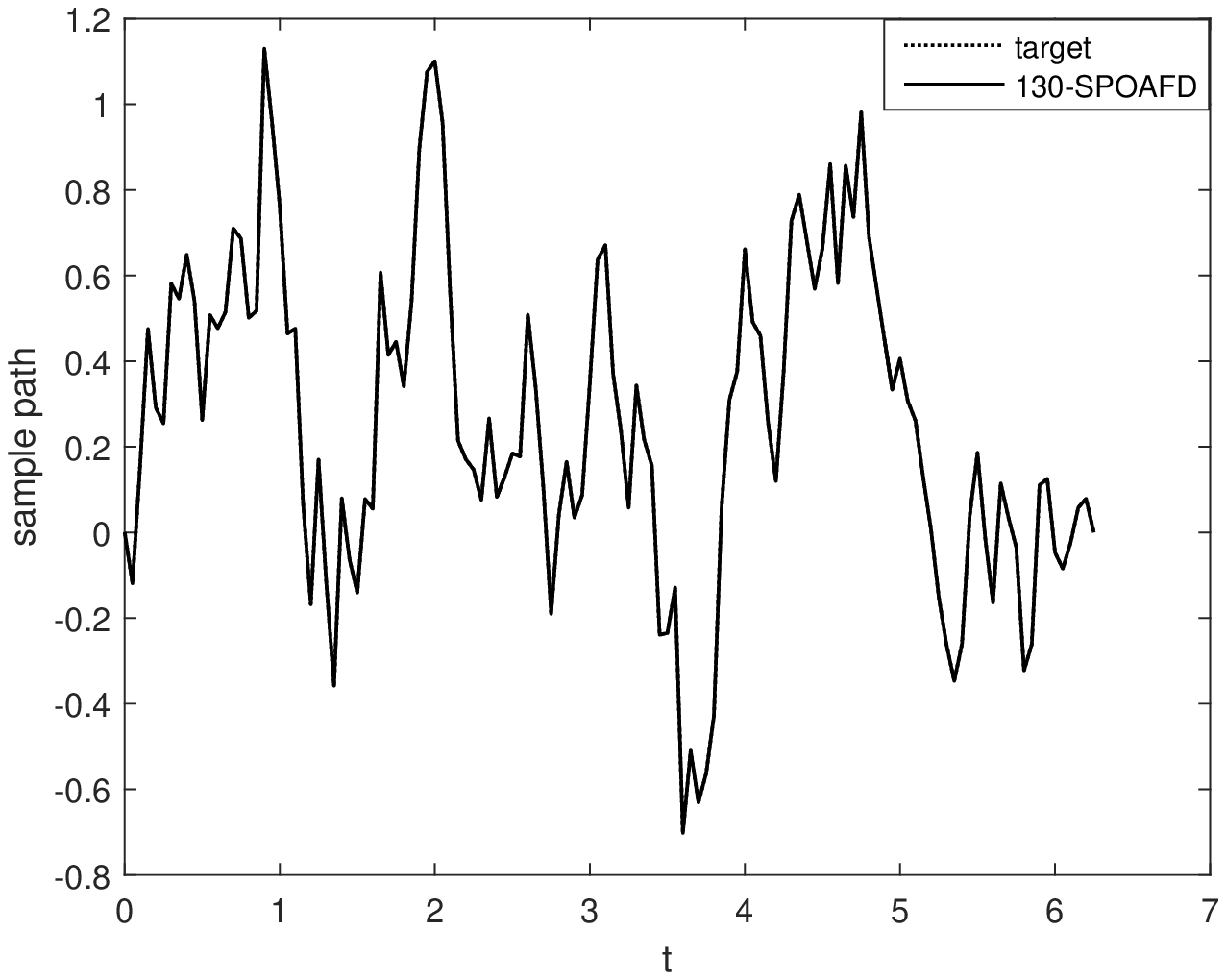}
		\centerline{{\scriptsize SPOAFD: 130 partial sum}}
	\end{minipage}
	\begin{minipage}[c]{0.3\textwidth}
		\centering
		\includegraphics[height=4.5cm,width=5cm]{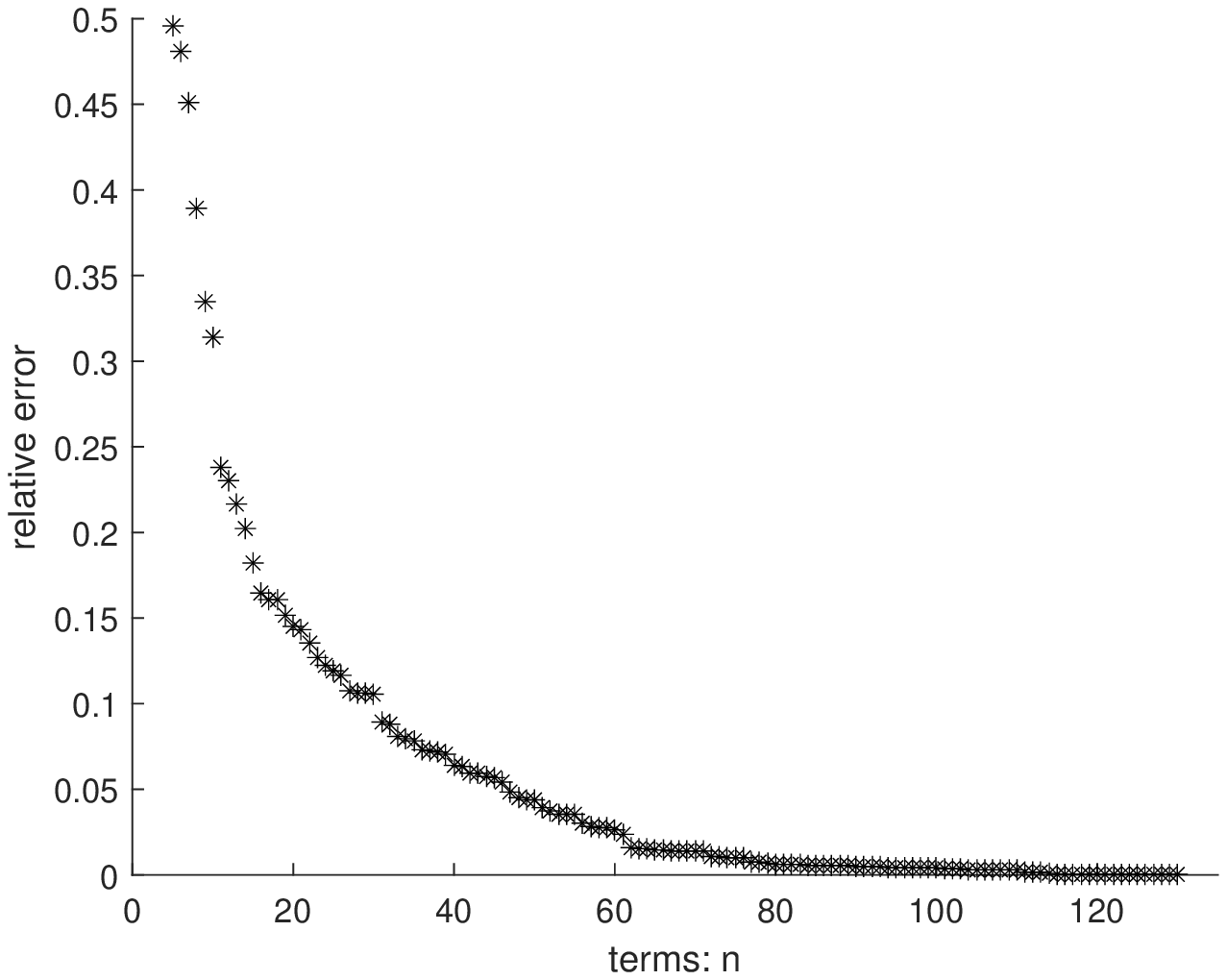}
		\centerline{{\scriptsize relative error}}
	\end{minipage}
	\caption{sample path II of Brownian bridge}
\label{figureB2}
\end{figure}

\begin{Exmp}\label{Exmp-Brownian}
In this example we decompose Brownian bridge.
We are considering the stochastic Dirichelet problem
\begin{equation}%\label{Equ-333}
\left\{
\begin{aligned}
&\Delta u(z,\omega)=0, ~ ~z \in B_1,\\
& u(e^{it})=W(t),~ ~ t\in [0,2\pi],
\end{aligned}
\right.
\end{equation}
where $B_1={\bf D}$ is the unit disc in the complex plane, and $W(t)$ is the Brownian bridge.
By using SPOAFD for the Poisson kernel in the unit disc we obtain an orthonormal system uniformly applicable for all Brownian bridges. We arbitrarily chose two sample pathes to test effectiveness of the obtained SPOAFD system. As shown in Figure \ref{figureB1} and \ref{figureB2}, the results are quite promising.
\end{Exmp}

  \section{Conclusion}

  Using the recently developed stochastic pre-orthogonal adaptive Fourier decomposition (SPOAFD) to the dictionaries of parameterised Poisson and
  heat kernels we obtain sparse representations of solutions of the Dirichlet boundary value and the Cauchy initial value problems with random data, including stochastic process cases. The numerical examples show promising convergence speed that implies fast convergence of the SPOAFD solutions series to the stochastic boundary and initial value problems. Theoretical and technical aspects of SPOAFD in comparison with existing major decomposition methods for stochastic processes will be studied in separate papers (\cite{KLvsAFD,OUprocess}).

\section*{Acknowledgments}

The study is supported by the Science and Technology Development Fund, Macau SAR (File no. 0123/2018/A3) and NSFC 12071437.

%\nolinenumbers

%\section*{References}


\begin{thebibliography}{10}
\bibitem{Axler_Harmonic-funct-thrm} S. Axler, P. Bourdon, W. Ramey, {\it Harmonic function theory}, Graduate Texts in Mathematics 137, 2013.

\bibitem{ACQS1} D. Alpay, F. Colombo, T. Qian, I. Sabadini, {\it Adaptive orthonormal systems for matrix-valued functions,} Proceedings of the American Mathematical Society, 2017, 145(5): 2089-2106.

\bibitem{ACQS2} D. Alpay,  F. Colombo, T. Qian, I. Sabadini, {\it Adaptive Decomposition: The Case of the Drury-Arveson Space,} Journal of Fourier Analysis and Applications, 2017, 23(6): 1426-1444.

\bibitem{Allen_1998} E. J. Allen, S. J. Novosel, Z. Zhang, {\it Finite element and difference approximation of some linear stochastic partial differential equations}, Stochastics-an International Journal of Probability and Stochastic Processes, 1998, 64(1): 117-142.

\bibitem{Babuska_1961} I. Babuska, {\it On randomized solutions of Laplace's equation}, asopis pro pestovani matematiky, 1961, 86(3).	

\bibitem{Babuka_Stochastic-Collocation2007} I. Babuka, F. Nobile, R. Tempone, {\it A Stochastic Collocation Method for Elliptic Partial Differential Equations with Random Input Data}, SIAM Journal on Numerical Analysis, 2007, 45(3): 1005-1034.	

\bibitem{Becus_rdm-Laplace-1977} G. A. Becus, {\it The random robin problem for laplace's equation}, Canadian Journal of Statistics, 2010, 5(1): 133-139.

\bibitem{Becus_rdm-heat-1977} G. A. Becus, {\it Random generalized solutions to the heat equation}, Journal of Mathematical Analysis and Applications, 1977, 60(1): 93-102.
	
\bibitem{Becus_approximations-rdm-heat-Equ-1978} G. A. Becus, {\it Solutions to the random heat equation by the method of successive approximations}, Journal of Mathematical Analysis and Applications, 1978, 64(2): 277-296.	
	
\bibitem{Becus_rdm-Laplace-1976} G. A. Becus, F. A. Cozzarelli, {\it The random steady state diffusion problem. I. Random generalized solutions to Laplace's equation}, SIAM J. Appl. Math, 1976, 31(1): 134-147.	

	
\bibitem{Boyaval_2009} S. Boyaval, C. L. Bris, Y. Maday, N. Nguyen, A. Patera, {\it A reduced basis approach for variational problems with stochastic parameters: Application to heat conduction with variable robin coefficient}, Computer Methods in Applied Mechanics and Engineering, 2009, 198(41-44): 3187-3206.
	
\bibitem{Cohen-2010}A. Cohen, R. DeVore, C. Schwab, {\it Convergence rates of best N-term Galerkin approximations for a class of elliptic SPDEs}, Foundations of Computational Mathematics, 2010, 10(6): 615-646.	
    	
\bibitem{Cameron-Martin1947} R. H. Cameron, W. T. Martin, {\it The orthogonal development of nonlinear functionals in a series of Fourier-Hermite functions}, Annals of Mathematics, 1947, 48(2): 385-392.

\bibitem{CP} R. Coifman, J. Peyri\'ere, {\it Phase unwinding, or invariant subspace decompositions of Hardy spaces,} Journal of Fourier Analysis and Applications, 2019, 25: 684-695.

\bibitem{CQT}  Q. H. Chen, T. Qian, L. H. Tan, {\it A Theory on Non-Constant Frequency Decompositions and Applications,} Advancements in Complex Analysis. Springer, Cham, 2020: 1-37.
\bibitem{CS} R. Coifman, S. Steinerberger, {\it Nonlinear phase unwinding of functions,} Journal of Fourier Analysis and Applications, 2017, 23: 778-809.

\bibitem{Choi_DO-schemes-2013} M. Choi, T. P. Sapsis, G. E. Karniadakis, {\it A convergence study for SPDEs using combined Polynomial Chaos and Dynamically-Orthogonal schemes}, Journal of Computational Physics, 2013, 245: 281-301.
			
\bibitem{Debussche_2007} A. Debussche, J. Printems, {\it Weak order for the discretization of the stochastic heat equation}, arXiv, 2007.	
	
\bibitem{Fefferman-Stein_Hardy-spc-1972} C. Fefferman, E. M. Stein, {\it $H^p$ Spaces of Several Variables}, Acta Mathematica, 1972, 129(1): 137-193.

\bibitem{Hausenblas_Galerkin-method2003} E. Hausenblas, {\it Approximation for Semilinear Stochastic Evolution Equations}, Potential Analysis, 2003, 18(2): 141-186.
	
\bibitem{Han-Lin_PDEs2011} Q. Han, F. H. Lin, {\it Elliptic partial differential equations}, Springer New York, 2011.
			
\bibitem{Hou_2006} T. Y. Hou, W. Luo, B. Rozovskii, H. M. Zhou, {\it Wiener chaos expansions and numerical solutions of randomly forced equations of fuid mechanics}, J Comput Phys, 2006, 216(2): 687-706.
	
\bibitem{Hytonen_AnalysisBanachSpaces_2016} T. Hytonen, J. V. Neerven, M. Veraar, L. Weis, {\it  Analysis in Banach Spaces, Volume I: Martingales and Littlewood-Paley Theory}, Springer, 2016.

\bibitem{Jardak_2002} M. Jardak, C. H. Su, G. E. Karniadakis, {\it Spectral polynomial chaos solutions of the stochastic advection equation}, Journal of Scientific Computing, 2002, 17(1): 319-338.	
\bibitem{Kaligotla_2012} S. Kaligotla, {\it Asymptotic Problems in Stochastic Partial Differential Equations: A Wiener Chaos Approach}, California: University of Southern California, 2012.
		
\bibitem{Kloeden-1992} P. E. Kloeden, E. Platen, {\it Numerical Solution of Stochastic Differential Equations}, Springer-Verlag Berlin Heidelberg, 1992.
	
\bibitem{Luo_2006} W. Luo, {\it Wiener chaos expansion and numerical solutions of stochastic partial differential equations}, Dissertations and Theses Gradworks, 2006.
		
\bibitem{LPS} G. J. Lord, C. E. Powell, T. Shardlow, {\it An Introduction to Computational Stochastic PDEs}, Cambridge University Press, 2014.	

\bibitem{Maruyama-Euler-Method-1955} G. Maruyama, {\it Continuous Markov processes and stochastic equations}, Rendiconti del Circolo Matematico di Palermo, 1955, 4(1): 48-90.
		
\bibitem{Machiels_Fourier-spectral-1998} L. Machiels, M. O. Deville, {\it Numerical Simulation of Randomly Forced Turbulent Flows}, Journal of Computational Physics, 1998, 145(1): 246-279.
	
\bibitem{Milstein-1995} G. N. Milstein, {\it Numerical Integration of Stochastic Differential Equations}, Kluwer Academic Publishers, 1995.

\bibitem{Milstein-2004} G. N. Milstein, M. V. Tretyakov, {\it Stochastic Numerics for Mathematical Physics}, Springer-Verlag Berlin Heidelberg, 2004.

\bibitem{Mugler_2013} A. Mugler, H. J. Starkloff, {\it On the convergence of the stochastic Galerkin method for random elliptic partial differential equations}, Esaim Mathematical Modelling and Numerical Analysis, 2013, 47(5): 1237-1263.
	
\bibitem{Nouy_spectral-decpst2007} A. Nouy, {\it A generalized spectral decomposition technique to solve a class of linear stochastic partial differential equations}, Computer Methods in Applied Mechanics and Engineering, 2007, 196(45-48): 4521-4537.
	
\bibitem{Nouy_spectral-decpst2008} A. Nouy, {\it Generalized spectral decomposition method for solving stochastic finite element equations: Invariant subspace problem and dedicated algorithms}, Computer Methods in Applied Mechanics and Engineering, 2008, 197(51-52): 4718-4736.

\bibitem{Nouy_Recent-developments-2009} A. Nouy, {\it Recent developments in spectral stochastic methods for the numerical solution of stochastic partial differential equations}, Archives of Computational Methods in Engineering, 2009, 16(3): 251-285.
		
\bibitem{Nobile-2008} F. Nobile, R. Tempone, C. Webster, {\it A sparse grid stochastic collocation method for partial differential equations with random input data}, SIAM J. Numer. Anal., 2008, 46(5): 2309-2345.
		
\bibitem{Park1998_K-L}H. M. Park, M. W. Lee, {\it An efficient method of solving the Navier-Stokes equation for the flow control}, Int. J. Numer. Methods Engrg, 1998, 41: 1131-1151.	

\bibitem{Park2000_K-L} H. M. Park, M. W. Lee, {\it Boundary control of the Navier-Stokes equation by empirical reduction of modes}, Comput. Methods Appl. Mech. Engrg, 2000, 188: 165-186.	

\bibitem{Qian2010} T. Qian,  {\it Intrinsic mono-component decomposition of functions: An advance of Fourier theory,} Mathematical Methods in the Applied Sciences, 2010, 33, 880-891.

\bibitem{Qianbook} T. Qian, {\it Adaptive Fourier Decomposition,}  Beijing: Science Press, 2015.

\bibitem{QianSAFD} T. Qian, {\it A Sparse Representation of Random Signals}, Mathematical Methods in the Applied Sciences, DOI: 10.1002/mma.8033, 2021.
	
\bibitem{Q2D} T. Qian, {\it Two-Dimensional Adaptive Fourier Decomposition}, Mathematical Methods in the Applied Sciences, 2016, 39(10) : 2431-2448.
	
\bibitem{Qian2018} T. Qian, {\it A novel Fourier theory on non-linear phase and applications,} Advances in Mathematics (China), 2018, {47}(3): 321-347.
\bibitem{KLvsAFD} T. Qian, et al.{\it Comparison of K-L and AFD type decompositions}, preprint.
\bibitem{Qu-Sps-represt-dirac} W. Qu, C. K. Chui, G. T. Deng, T. Qian, {\it sparse representation of approximation to identity}, Analysis and Applications, 2021.

\bibitem{qu2018} W. Qu, P. Dang, {\it Rational approximation in a class of weighted Hardy spaces,} Complex Analysis and Operator Theory, 2019, 13(4): 1827-1852.
	
\bibitem{qu2019} W. Qu, P. Dang, {\it Reproducing kernel approximation in weighted Bergman spaces: Algorithm and applications}, Mathematical Methods in the Applied Sciences, 2019, 42(12): 4292-4304.

\bibitem{QQ} W. Qu, T. Qian, H. C. Li, K. H. Zhu, {\it Best kernel approximation in Bergman spaces,} Applied Mathematics and Computation, 2022, 416:126749.

\bibitem{QTa} T. Qian, L. H. Tan, J. C. Chen, {\it A class of iterative greedy algorithms related to Blaschke product}, Science China Mathematics, 2021, https://doi.org/10.1007/s11425-020-1706-5.

\bibitem{QWM} T. Qian, J. Z. Wang, W. X. Mai, {\it An Enhancement Algorithm for Cyclic Adaptive Fourier Decomposition}, Applied and Computational Harmonic Analysis, 2019, 47(2): 516-525.
	
\bibitem{QSW} T. Qian, W. Sproessig, J. X. Wang, {\it Adaptive Fourier decomposition of functions in quaternionic Hardy spaces, Mathematical Methods in the Applied Sciences}, 2012, 35(1): 43-64.

\bibitem{QWa} T.~Qian, Y. B. Wang, {\it Adaptive Fourier series-a variation of greedy algorithm}, Advances in Computational Mathematics, 2011, 34~(3): 279-293.

\bibitem{Sabelfeld2008-expansion-random-boundary-PDE} K. Sabelfeld, {\it Expansion of random boundary excitations for elliptic PDEs}, Monte Carlo Methods and Applications, 2008, 13(5-6): 405-453.	
	
\bibitem{Sapsis_DO-2009} T. P. Sapsis, P. Lermusiaux, {\it Dynamically orthogonal field equations for continuous stochastic dynamical systems}, Physica D Nonlinear Phenomena, 2009, 238(23-24): 2347-2360.	
	
\bibitem{Schwab-2011} C. Schwab, C. Gittelson, {\it Sparse tensor discretizations of high-dimensional parametric and stochastic PDEs}, Acta Numerica, 2011, 20(20): 291-467.

\bibitem{Schwab2003_sparse-FE-PDE-stochastic} C. Schwab, R. A. Todor, {\it Sparse finite elements for elliptic problems with stochastic loading}, Numerische Mathematik, 2003, 95(4): 707-734.

\bibitem{Saitoh_RKHS-2016} S. Saitoh, Y. Sawano, {\it Theory of Reproducing Kernels and Applications}, Developments in Mathematics, 2016, 44.	
	
\bibitem{Stein-Weiss1971} E. M. Stein, G. Weiss, {\it Introduction to Fourier Analysis on Euclidean Space}, Princeton University Press, 1971.

	
\bibitem{Tasaka_FE-method-Rdm-heat-Equ-1983} S. Tasaka, {\it Convergence of statistical finite element solutions of the heat equation with a random initial condition}, Computer Methods in Applied Mechanics and Engineering, 1983, 39(2): 131-136.

\bibitem{TY} V. N. Temlyakov, {\it Greedy approximation}, Cambridge University Press, 2011.

\bibitem{Todor-2007} R. A. Todor, C. Schwab, {\it Convergence rates for sparse chaos approximations of elliptic problems with stochastic coefficients}, IMA J. Numer. Anal.,2007, 27(2): 232-261.
	    	
\bibitem{Walsh_Finite-Element-2005} J. B. Walsh, {\it Finite Element Methods for Parabolic Stochastic PDE's}, Potential Analysis, 2005, 23(1): 1-43.

\bibitem{WQ2020} Y. B. Wang, T. Qian, {\it Pseudo-hyperbolic distance and $n$-best rational approximation in $H^2$ space}, Mathematical Methods in the Applied Sciences, 2021, 44(11): 8497-8504.

\bibitem{Xiu2010-numerical-methods} D. B. Xiu, {\it  Numerical Methods for Stochastic Computations}, Princeton University Press, 2010.

\bibitem{Xiu_2002} D. B. Xiu, G. E. Karniadakis, {\it  Modeling uncertainty in steady state diffusion problems via generalized polynomial chaos}, Computer Methods in Applied Mechanics and Engineering, 2002, 191(43): 4927-4948.	
	
\bibitem{Yosida1978_FunctionalAnalysis} B. Yosida, {\it Functional Analysis}, Springer-Verlag, 1978.		
	
\bibitem{Yannacopoulos_2011} A. N. Yannacopoulos, N. E. Frangos, I. Karatzas, {\it Wiener chaos solutions for linear backward stochastic evolution equations}, SIAM J Math Anal, 2011, 43(1): 68-113.	

\bibitem{Georgi_2017} Z. Q. Zhang, G. E. Karniadakis, {\it  Numerical Methods for Stochastic Partial Differential Equations with White Noise}, Springer, 196, 2017.

\bibitem{OUprocess}	Y. Zhang, {\it Sparse repersetations of Ornstein-Uhlenbeck process}, preprint.
\end{thebibliography}
\end{document}